\documentclass[preprint,12pt]{elsarticle}
\usepackage{graphicx}
\usepackage{amssymb}
\usepackage{amsthm}
\usepackage{graphics}
\usepackage{epsfig}
\usepackage{amssymb}
\usepackage{hyperref}
\usepackage{graphicx}
\usepackage{subfigure}
\usepackage{float}
\usepackage{blkarray}
\usepackage{blkarray,bigstrut}
\usepackage{color}
\usepackage{tikz-network}
\usepackage{tikz}
\usepackage{tikz,pgfplots}
\pgfplotsset{compat=1.18}
\usetikzlibrary{patterns}
\usepackage[ruled,vlined]{algorithm2e}
\usepackage{caption}
\numberwithin{equation}{section}
\DeclareCaptionLabelFormat{cont}{#1~#2\alph{ContinuedFloat}}
\usepackage{fullpage}
\usepackage{amsmath}
\usepackage{amsfonts}
\usepackage{enumerate}
\usepackage{enumitem} 
\usepackage{tikz}
\usepackage{array}
\usepackage{bigdelim}
\usepackage{nicematrix}
\usepackage{mathtools}
\usepackage{blkarray, bigstrut}
\usepackage{gauss}
\usepackage{xcolor}
\usepackage{tabularray}
\usepackage{multicol}
\usepackage{multirow}
\usepackage{makecell}
\usepackage[most]{tcolorbox}
\UseTblrLibrary{booktabs}

\newcommand{\pf}{\noindent {\bf Proof: }}

\newtheorem*{theoremaux}{Theorem \theoremauxnum}
\gdef\theoremauxnum{1}

\newenvironment{theoremff}[2][]{%
	\def\theoremauxnum{\ref{#2}}
	\begin{theoremaux}[#1]
	}{%
\end{theoremaux}
}

\newtheorem{lemma}{\bf Lemma}[section]
\newtheorem{theorem}{\bf Theorem}[section]

\newtheorem{proposition}[lemma]{\bf Proposition}
\newtheorem{corollary}[lemma]{\bf Corollary}
\newtheorem{definition}{\bf Definition}[section]
\theoremstyle{definition}
\newtheorem{remark}{\bf Remark}[section]
\newtheorem{example}{\bf Example}[section]

\journal{~}

\begin{document}

\begin{frontmatter}
\title{{Integral quasi-strongly regular Cayley graphs}}

\author[1]{Sauvik Poddar}
\ead{sauvikpoddar1997@gmail.com}

\author[2]{Sucharita Biswas\corref{cor1}}
\ead{biswas.sucharita56@gmail.com}

\author[1]{Angsuman Das}
\ead{angsuman.maths@presiuniv.ac.in}

\address[1]{Department of Mathematics, Presidency University, Kolkata 700073, India\\
}
\address[2]{Department of Mathematics, Indian Institute of Technology Bombay, Mumbai 400076, India}
\cortext[cor1]{Corresponding author}

\begin{abstract}
Quasi-strongly regular graphs form a significant generalization of strongly regular graphs. We study the eigenvalues of a family of such graphs, $\Gamma_H(G)$, constructed from a finite group $G$ and a subgroup $H$. Our main results include a sufficient condition for $\Gamma_H(G)$ to be integral and an explicit computation of its entire spectrum when $H$ is normal, revealing that the spectrum in this case depends only on $|G|$ and the index $[G:H]$.
\end{abstract}


\end{frontmatter}

\section{Introduction}

Strongly regular graphs (SRGs), introduced by Bose \cite{bose1963strongly}, are a cornerstone of algebraic and spectral graph theory, with a rich and extensive literature. This has naturally led to several significant generalizations, such as \textit{Deza graphs} \cite{deza1994ridge,erickson1999deza}, \textit{generalized strongly regular graphs} \cite{huo2017subconstituents} and \textit{quasi-strongly regular graphs} \cite{golightly1997family,xie2021quasi,xie2023quasi}. This paper focuses on the spectral properties of a specific family of quasi-strongly regular Cayley graphs. 

\begin{definition}
Let $G$ be a finite group and $S\subseteq G\setminus\lbrace{1}\rbrace$ be an inverse-closed subset, i.e., $S=S^{-1}$, where $S^{-1}\coloneqq\lbrace{s^{-1}\mid s\in S}\rbrace$. A Cayley graph $\operatorname{Cay}(G,S)$ of G with respect to $S$ is a graph having $G$ as its vertex set, where two distinct vertices $u,v\in G$ are adjacent if $uv^{-1}\in S$.
\end{definition}

The set $S$ is said to be the \textit{connection set} for
$\operatorname{Cay}(G,S)$ and the size of $S$ is called the \textit{valency} of the Cayley graph. A set $S\subseteq G$ is said to be \textit{normal} in $G$ if $gSg^{-1}=S$ for all $g\in G$. A Cayley graph is said to be \textit{normal} if the connection set is normal.

The adjacency spectrum of a graph $\Gamma$ is the multiset of its adjacency eigenvalues and is denoted by
$$\operatorname{Spec}(\Gamma)=\begin{pmatrix}
\lambda_1 & \lambda_2 & \cdots & \lambda_k\\
m_1 & m_2 & \cdots & m_k
\end{pmatrix},$$
where the first row denotes the distinct adjacency eigenvalues of $\Gamma$ and the second row denotes their respective multiplicities. It can be expressed alternatively as
$$\operatorname{Spec}(\Gamma)=\left[[\lambda_1]^{m_1},[\lambda_2]^{m_2},\ldots,[\lambda_k]^{m_k}\right],$$
where the exponent denotes the multiplicities of each eigenvalue. We usually omit $m_i$ if $m_i=1$ for some $i$.

Two graphs $\Gamma_1$ and $\Gamma_2$ are said to be \textit{isospectral} (or \textit{cospectral}) if $\operatorname{Spec}(\Gamma_1)=\operatorname{Spec}(\Gamma_2)$ as multisets. An eigenvalue of $\Gamma$ is said to be \textit{simple}, if its multiplicity is $1$. A graph $\Gamma$ is said to be \textit{singular} if $0\in \operatorname{Spec}(\Gamma)$. The study of non-isomorphic, isospectral graph pairs occupies a prominent place in spectral graph theory research. Our results, including Theorem \ref{qsrg-isospectral}, give an efficient and elegant method to generate many such pairs.

A graph $\Gamma$ is called \textit{integral} if $\operatorname{Spec}(\Gamma)$ consists entirely of integers. The quest of characterizing integral graphs dates back to $1973$ when Harary and Schwenk \cite{harary1974graphs} proposed the problem of classifying all integral graphs. This problem, in spite of its simple and straightforward description, turns out to be exceedingly difficult.

In the literature of spectral graph theory, one of the main research areas concerns classifying integral graphs, in particular Cayley graphs. Although Alperin and Peterson \cite{alperin2012integral} have completely characterized the integrality condition for Cayley graphs over abelian groups, and later works extended these results to specific non-abelian families \cite{cheng2019integral,cheng2023integral,huang2021integral,lu2018integral}, a general understanding remains elusive. A key recent result by Godsil and Spiga \cite{godsil2025integral} provides a necessary and sufficient condition for the integrality of normal Cayley graphs. In this paper, we establish a sufficient condition for the integrality of a family of Cayley graphs that are not necessarily normal.

\subsection{Some definitions and results}

\medskip

\begin{definition}\cite{goldberg2006quasi}
A strongly regular graph with parameters $(n,k,a,c)$, denoted by \\$SRG(n,k,a,c)$, is a $``k"$ regular graph on $``n"$ vertices such that any two adjacent vertices have $``a"$ many common neighbours and any non-adjacent vertices have $``c"$ many common neighbours.
\end{definition}

\begin{definition}\cite{golightly1997family}
A quasi-strongly regular graph with parameters $(n,k,a;c_1,c_2,\ldots,c_p)$, denoted by $QSRG(n,k,a;c_1,c_2,\ldots,c_p)$, is a $``k"$ regular graph on $``n"$ vertices such that any two adjacent vertices have $``a"$ many common neighbours and any two non-adjacent vertices have $``c_i"$ many common neighbours for some $1\le i\le p$, where $c_i$’s are distinct non-negative integers. The set of parameters $\lbrace{c_1,c_2,\ldots,c_p}\rbrace$ is called the $c$-set of
$QSRG(n,k,a;c_1,c_2,\ldots,c_p)$ and $p$ is called the grade.
\end{definition}

A QSRG of grade $1$ is precisely an SRG. A QSRG is called proper if its grade $p>1$. In \cite{biswas2025family}, the authors constructed a family of quasi-strongly regular Cayley graphs by generalizing a well-known construction of SRGs from finite groups and characterized its full automorphism group and other transitivity properties.

\begin{definition}\cite{biswas2025family}
Let $G$ be a finite group of order $n\ge 5$, $H$ be a proper subgroup of $G$ and $S_H=\lbrace{(g,1),(1,g),(g,g)\mid g\in G\setminus H}\rbrace$. Then $\Gamma_{H}(G)$ is defined to be the Cayley graph of $G\times G$ with respect to the connection set $S_H$, i.e., $\Gamma_{H}(G)\coloneqq\operatorname{Cay}(G\times G,S_H)$.
\end{definition}

If $H=G$, then $\Gamma_{H}(G)$ becomes an edgeless graph and if $H=\lbrace{1}\rbrace$, then $\Gamma_{H}(G)$ belongs to a well-known family of strongly regular graphs with the parameters $(n^2,3n-3,n,6)$ as shown in \cite[Page 27, Theorem 3.22]{nica2018brief}, which yields an integral family of Cayley graphs. If $H$ is a proper, non-trivial subgroup of $G$, then $\Gamma_{H}(G)$ is not strongly regular. The vertices in $S_H$ can be regarded as Type $1$, Type $2$ and Type $3$ for the elements $(g,1),(1,g)$ and $(g,g)$, respectively, where $g\in G\setminus H$. Note that $\Gamma_{H}(G)$ is not a normal Cayley graph unless $G$ is abelian. The authors in \cite{biswas2025family} showed that $\Gamma_H(G)$ is a quasi-strongly regular graph of grade $3,4$ and $5$.

\begin{theorem}\cite[Theorem 2.1]{biswas2025family}
For any subgroup $H$ of a finite group $G$, any two adjacent vertices of $\Gamma_H(G)$ have $|G|-2|H|+2$ many common neighbours, that is, $a=|G|-2|H|+2$.
\end{theorem}

\begin{theorem}\cite[Theorem 2.2]{biswas2025family}
Let $H$ be a non-trivial proper subgroup of $G$ and $[G:H]=\ell$.
\begin{enumerate}
\item If $\ell=2$, then $c$-set $=\lbrace{0,2,|G|-|H|}\rbrace$.
\item If $\ell>2$ and $H$ is normal in $G$, then
\begin{itemize}
\item $|H|=2$ implies $c$-set $=\lbrace{2,6,|G|-|H|}\rbrace$.
\item $|H|>2$ implies $c$-set $=\lbrace{0,2,6,|G|-|H|}\rbrace$.
\end{itemize}
\item If $\ell>2$ and $H$ is not normal in $G$, then
\begin{itemize}
\item $|H|=2$ implies $c$-set $=\lbrace{2,4,6,|G|-|H|}\rbrace$.
\item $|H|>2$ implies $c$-set $=\lbrace{0,2,4,6,|G|-|H|}\rbrace$.
\end{itemize}
\end{enumerate}
\end{theorem}

Another aspect of this article concerns the spectrum of $\Gamma_H(G)$, in light of the quasi-strong regularity. The spectrum of a strongly regular graph is well known to consist of exactly three distinct eigenvalues; conversely, any connected regular graph with exactly three distinct eigenvalues is strongly regular. In contrast, very little is known about the spectra of quasi-strongly regular graphs of grade greater than $2$, as determining their eigenvalues in terms of parameters is significantly more difficult. Some spectral results exist for QSRGs of grade $2$ (e.g., \cite{guo2022some,xie2024spectra}), but to the best of our knowledge, the spectral properties of QSRGs of higher grade have not been studied extensively. This paper aims to address this gap by deriving the spectrum of the family $\Gamma_H(G)$ with grade greater than $2$. 

\subsection{Our contribution}

\medskip

Our contribution in this paper is mainly two-fold: we prove our two main results in Sections \ref{section-3} and \ref{section-5}. Section \ref{section-3} provides a sufficient condition for the integrality of the graph $\Gamma_H(G)$. In particular, we prove the following:

\begin{theoremff}{H-normal-graph-integral}
Let $G$ be a finite group and $H\lhd G$. Then $\Gamma_H(G)$ is integral.
\end{theoremff}

In Section \ref{section-5}, we determine the complete spectrum of the graph $\Gamma_H(G)$, when $H$ is a normal subgroup of $G$, in terms of $|G|,|H|$ and $[G:H]$, and thereby prove the following:

\begin{theoremff}{qsrg-isospectral}
Let $G_1$ and $G_2$ be two groups of same order. Let $H_1\lhd G_1$ and $H_2\lhd G_2$ be such that $|H_1|=|H_2|$. Then the graphs $\Gamma_{H_1}(G_1)$ and $\Gamma_{H_2}(G_2)$ are isospectral.
\end{theoremff}

In other words, for a group $G$ and its normal subgroup $H$, the spectrum of $\Gamma_H(G)$ is independent of the group structures of $G$ and $H$, and solely depends upon the orders of $G$ and $H$. 

Additionally, in Section \ref{section-4}, we obtain some general results concerning the spectrum of $\Gamma_H(G)$.

\section{Definitions and Preliminaries}

We consider only simple (i.e., without loops and multiedges), undirected, finite graphs. All groups considered here are finite. To denote sets and multisets, we use curly brackets $\lbrace{}\rbrace$ and square brackets $[~]$, respectively. The interpretation of the symbol $\le$ varies depending on the context of the objects involved. Specifically, if $V$ and $W$ are considered as vector spaces, then $W\le V$ will denote that $W$ is a subspace of $V$, while in the context of groups, $H\le G$ will denote that $H$ is a subgroup of $G$. For numbers, $\le$ and $\ge$ will bear their usual meanings respectively. Since the cases $H=G$ and $H=\lbrace{1}\rbrace$ are trivial, we shall, without loss of generality, consistently use the notation $H\le G$ to denote that $H$ is a proper, non-trivial subgroup of $G$, while considering $\Gamma_H(G)$. If $\lambda$ is an eigenvalue of the graph $\Gamma$, we denote the multiplicity of $\lambda$ in $\Gamma$ by $m(\lambda)$. We denote by $\mathcal{K}(\Gamma)$ and $\kappa(\Gamma)$, the set and the number of distinct eigenvalues of the graph $\Gamma$, respectively.

\subsection{Basics of representations and characters}

\medskip

We now recall some basics of representation and character theory of a finite group and also review some important results about the spectrum of Cayley graphs. For more information about the representation theory and character theory of finite groups, one can refer to \cite{isaacs1994character,serre1977linear,steinberg2012representation}.

Let $G$ be a finite group and $V$ be a finite-dimensional vector space over $\mathbb{C}$. A \textit{representation} $(\rho,V)$ of $G$ is a group homomorphism $\rho:G\to GL(V)$, where $GL(V)$ denotes the group of all invertible linear transformations of $V$. We can simply denote a representation as $\rho$ if $V$ is already understood from the context. The \textit{degree} of $\rho$, denoted by $\deg(\rho)$, is $\dim(V)$, the dimension of $V$. Two representations $(\rho,V)$ and $(\varphi,W)$ of $G$ are called \textit{equivalent}, written as $\rho\sim\varphi$, if there exists an isomorphism $T:V\to W$ such that $T\rho(g)=\varphi(g)T$ for all $g\in G$. A \textit{morphism} from $\rho$ to $\varphi$ is a linear map $T:V\to W$ such that $T\rho(g)=\varphi(g)T$ for all $g\in G$. The set of all morphisms from $\rho$ to $\varphi$ is denoted by $\operatorname{Hom}_G(\rho,\varphi)$ and is a subspace of $\operatorname{Hom}(V,W)$, the set of all linear transformations from $V$ to $W$. If $T\in\operatorname{Hom}_G(\rho,\varphi)$ is invertible, then $\rho\sim\varphi$ and $T$ is an equivalence (or isomorphism). 

The character $\chi_\rho:G\to\mathbb{C}$ of $\rho$ is a map defined as 
$$g\mapsto\operatorname{Tr}(\rho(g)),~\text{for $g\in G$},$$
where $\operatorname{Tr}(\rho(g))$ is the trace of the representation matrix of $\rho(g)$ with respect to some basis of $V$. The \textit{degree} of the character $\chi_\rho$, denoted $\deg(\chi_{\rho})$, is the degree of $\rho$ and is equal to $\chi_\rho(1)$. A subspace $W\le V$ is said to be \textit{$G$-invariant} if $\rho(g)w\in W$ for each $g\in G$ and $w\in W$. If $W$ is a $G$-invariant subspace of $V$, then the restriction of $\rho$ on $W$, i.e., $\rho\big|_W: G \to GL(W)$, is a representation of $G$ on $W$. The subspaces $\{0\}$ and $V$ are always $G$-invariant, and are referred to as trivial subspaces. We say that $(\rho,V)$ is an \textit{irreducible representation} and $\chi_\rho$ an \textit{irreducible character} of $G$, if $V$ has no non-trivial $G$-invariant subspaces. We denote by $\operatorname{IRR}(G)$ and $\operatorname{Irr}(G)$ the complete set of inequivalent (complex) irreducible representations of $G$ and the complete set of inequivalent irreducible (complex) characters of $G$, respectively.

\subsection{Group algebra and regular representation}

\medskip

The group algebra $\mathbb{C}[G]$ is a vector space whose basis consists of the elements of $G$, i.e.,
$$\mathbb{C}[G]=\left\{\sum_{g\in G}c_gg\mid c_g\in\mathbb{C}\right\}.$$
The \textit{(left) regular representation} of $G$ is the homomorphism $L:G \to GL(\mathbb{C}[G])$ defined by
$$L(g)\left(\sum_{h\in G}c_hh\right)=\sum_{h\in G}c_hgh=\sum_{x\in G}c_{g^{-1}x}x,$$
for each $g\in G$. The character $\chi_L$ of the regular representation is given by
$$\chi_L(g)=\begin{cases}
    |G|,& \text{if $g=1$.}\\
    0, & \text{if $g\ne 1$.}
\end{cases}$$

\begin{lemma}\label{lem:cayley1}\cite{serre1977linear}
Let $G$ be a finite group and $\operatorname{IRR}(G)=\lbrace{\rho_1,\dots,\rho_h}\rbrace$. Let $L$ be the regular representation of $G$. Then
$$L\sim \bigoplus_{i=1}^{h}\rho_i^{\oplus d_i},$$
where $d_i=\deg(\rho_i)$, $1\le i\le h$ and
$$\rho_i^{\oplus d_i}=\underbrace{\rho_i\oplus\cdots\oplus\rho_i}_{d_i~\text{times}}.$$
\end{lemma}

Let $R(g)$ denote the representation matrix corresponding to $L(g)$ for $g \in G$. Babai \cite{babai1979spectra} showed that the adjacency matrix of $\operatorname{Cay}(G,S)$ can be expressed in terms of $R(g)$.

\begin{lemma}\label{lemma-babai}\cite{babai1979spectra}
Let $G$ be a finite group of order $n$, and let $S \subseteq G\setminus\lbrace{1}\rbrace$ be such that $S = S^{-1}$. Then the adjacency matrix $A$ of $\operatorname{Cay}(G,S)$ can be expressed as $A = \sum_{s \in S} R(s)$, where $R(s)$ is the representation matrix corresponding to $L(s)$.
\end{lemma}

We denote by $R_i(g)$, the representation matrix corresponding to $\rho_i(g)$ for $g \in G$ and $1 \leq i \leq h$. Then by Lemma \ref{lem:cayley1}, there exists an orthogonal matrix $P$ such that
\begin{equation}\label{orthogonal-matrix-similar}
PR(g)P^{-1}=\bigoplus_{i=1}^{h} R_i(g)^{\oplus d_i}, 
\end{equation}
for each $g \in G$. Therefore, we have
\begin{equation}{\label{eqn:unique}}
PAP^{-1} = P\left(\sum_{s\in S}R(s)\right)P^{-1}=\bigoplus_{i=1}^{h}\left(\sum_{s \in S}R_i(s)\right)^{\oplus d_i}.
\end{equation}
Let $\rho\in\operatorname{IRR}(G)$ and $\chi\in\operatorname{Irr}(G)$. For $S\subseteq G$, define the notation $\rho(S)\coloneqq\sum_{s\in S}\rho(s)$ and $\chi(S)\coloneqq\sum_{s\in S}\chi(s)$. 

\begin{lemma}\label{chi-afford-rho-lemma}
Let $G$ be a finite group. Let $\rho\in\operatorname{IRR}(G)$ and $\chi\in\operatorname{Irr}(G)$ be the corresponding irreducible character. Let $S\subseteq G$ be normal in $G$. Then 
$$\rho(S)=\frac{\chi(S)}{\chi(1)}I_{\chi(1)}.$$
\end{lemma}
\pf This is an easy consequence of \cite[Lemma $5$]{diaconis1981generating}, which involves Schur's lemma.\qed 

\subsection{Tensor product of representations, fixed subspace and dual representation}

\medskip

Let $G_1$ and $G_2$ be two groups and $(\rho_1,V_1)$, $(\rho_2,V_2)$ be two representations of $G_1$ and $G_2$ respectively. The \textit{tensor product} of the representations $(\rho_1,V_1)$ and $(\rho_2,V_2)$ is the representation $(\rho_1\otimes\rho_2,V_1\otimes V_2)$ of $G_1\times G_2$ defined as
$$\rho_1\otimes\rho_2:G_1\times G_2\to GL(V_1\otimes V_2)$$
$$(\rho_1\otimes\rho_2)(g_1,g_2)=\rho_1(g_1)\otimes\rho_2(g_2).$$
If $\chi_1$ and $\chi_2$ are the characters corresponding to the representations $\rho_1$ and $\rho_2$, the character of the representation $\rho_1\otimes\rho_2$ becomes $\chi_1\chi_2$, where 
$$(\chi_1\chi_2)(g_1,g_2)=\chi_1(g_1)\chi_2(g_2).$$
The irreducible representations and irreducible characters of $G_1\times G_2$ are respectively given by 
$$\operatorname{IRR}(G_1\times G_2)=\lbrace{\rho\otimes\psi\mid \rho\in\operatorname{IRR}(G_1),\psi\in\operatorname{IRR}(G_2)}\rbrace,$$
$$\operatorname{Irr}(G_1\times G_2)=\lbrace{\chi\phi\mid \chi\in\operatorname{Irr}(G_1),\phi\in\operatorname{Irr}(G_2)}\rbrace.$$
Let $\chi,\psi\in\operatorname{Irr}(G)$. The inner product of characters is defined by
$$\langle{\chi,\psi}\rangle=\frac{1}{|G|}\sum_{g\in G}\chi(g)\overline{\psi(g)}.$$
Let $H\le G$ and $(\rho,V)$ be a representation of $G$. The \textit{$G$-fixed} and \textit{$H$-fixed} subspaces of $V$ are respectively defined as
$$V^G\coloneqq\lbrace{v\in V\mid \rho(g)v=v, \forall g\in G}\rbrace,$$
$$V^H\coloneqq\lbrace{v\in V\mid \rho(h)v=v, \forall h\in H}\rbrace.$$
It can be shown that $V^G$ is a $G$-invariant subspace of $V$. Clearly, for any subgroup $H\le G$, $V^G\le V^H\le V$. The dimension of $V^G$ is the rank of the operator $P_G:V\to V$ where
$$P_G\coloneqq\frac{1}{|G|}\sum_{g\in G}\rho(g).$$
The image of $P_G$ is $V^G$. Also, it follows immediately that $P_G^2=P_G$, i.e., $P_G$ is a projection operator. Then we have $\dim\left(V^G\right)=\operatorname{Rank}(P_G)$. Since the rank of the projection operator is equal to its trace, we further have 
$$\operatorname{Rank}(P_G)=\operatorname{Tr}(P_G)=\frac{1}{|G|}\sum_{g\in G}\operatorname{Tr}(\rho(g)).$$
Hence the dimension of $V^G$ can be obtained in terms of character,
$$\dim\left(V^G\right)=\frac{1}{|G|}\sum_{g\in G}\chi_{\rho}(g).$$
Let $(\rho,V)$ be a representation of $G$. The \textit{dual} $\rho^*$ of the representation $\rho$ is defined as,
$$\rho^*:G\to GL(V^*)$$
$$\rho^*(g)=\rho(g^{-1})^T,$$
where $V^*=\operatorname{Hom}(V,\mathbb{C})$ is the dual space of $V$.

A representation is irreducible if and only if its dual is irreducible. Also for any representation $\rho$, $(\rho^*)^*=\rho$ and $\deg(\rho)=\deg(\rho^*)$. If $\rho$ is a representation with character $\chi_{\rho}$, then $\chi_{\rho^*}=\overline{\chi_{\rho}}$. Let $(\rho,V)$ and $(\varphi,W)$ be two irreducible representations of $G$. Then $\dim(V\otimes W)^G=\dim\operatorname{Hom}_G(\rho^*,\varphi)$. It is known by Schur's lemma that
$$\dim\operatorname{Hom}_G(\rho^*,\varphi)=\begin{cases}
    1,&\text{if $\rho^*\sim\varphi$.}\\
    0,&\text{if $\rho^*\nsim\varphi$.}
\end{cases}$$
Alternatively, in terms of characters, this is equivalent to
$$\dim(V\otimes W)^G=\frac{1}{|G|}\sum_{g\in G}\chi_{\rho}(g)\chi_{\varphi}(g)=\langle{\chi_{\rho},\overline{\chi_{\varphi}}}\rangle=\langle{\chi_{\rho},\chi_{\varphi^*}}\rangle=\begin{cases}
    1,&\text{if $\rho\sim\varphi^*$.}\\
    0,&\text{if $\rho\nsim\varphi^*$.}
\end{cases}$$
Let $(\rho,V)$ be a representation of $G$ and $\chi$ be its corresponding character. The \textit{kernel} of a representation $\rho$, defined as 
$$\ker(\rho)\coloneqq\lbrace{g\in G\mid \rho(g)=I}\rbrace$$
is a normal subgroup of $G$.  Since $g\in\ker(\rho)$ if and only if $\chi(g)=\chi(1)$, the kernel of a character $\chi$ can be likewise defined as 
$$\ker(\chi)\coloneqq\lbrace{g\in G\mid \chi(g)=\chi(1)}\rbrace.$$
Let $H\le G$. We say that a representation $\rho$ (similarly, a character $\chi$) is \textit{trivial} on $H$ if $H\subseteq\ker(\rho)$ (similarly, $H\subseteq\ker(\chi)$). If $H\lhd G$, it can be shown that the irreducible characters of $G/H$ are precisely the irreducible characters of $G$ trivial on $H$, i.e., 
$$\operatorname{Irr}(G/H)=\lbrace{\chi\in\operatorname{Irr}(G)\mid H\subseteq\ker(\chi)}\rbrace.$$
Similarly, in the notion of representation, we have 
$$\operatorname{IRR}(G/H)=\lbrace{\rho\in\operatorname{IRR}(G)\mid H\subseteq\ker(\rho)}\rbrace.$$

\begin{lemma}\label{closed-under-dual-lemma}
Let $H\lhd G$ and $\rho\in\operatorname{IRR}(G)$. Then $\rho\in\operatorname{IRR}(G/H)$ if and only if $\rho^*\in\operatorname{IRR}(G/H)$. 
\end{lemma}
\pf The result follows immediately from the fact that $H\subseteq\ker(\rho)$ if and only if $H\subseteq\ker(\rho^*)$.\qed

\begin{lemma}\label{G-H-char-lemma}
Let $G$ be a finite group and $H\lhd G$. Then for every $\chi\in\operatorname{Irr}(G)$,
\begin{enumerate}[label=(\roman*)]
\item $$\chi(G)\coloneqq\sum_{g\in G}\chi(g)=\begin{cases}
    |G|, & \text{if $\chi={\bf{1}}_G$.}\\
    0, & \text{if $\chi\ne{\bf{1}}_G$.}
\end{cases}$$\\
\item $$\chi(H)\coloneqq\sum_{h\in H}\chi(h)=\begin{cases}
    |H|\chi(1), & \text{if $H\subseteq\ker(\chi)$.}\\
    0, & \text{if $H\nsubseteq\ker(\chi)$.}
\end{cases}$$
\end{enumerate}
\end{lemma}
\pf For (i), note that $\chi(G)=|G|\langle{\chi,{\bf{1}}_G}\rangle$, where ${\bf{1}}_G$ is the trivial character of $G$. The result follows from the orthogonality relations of characters.

\medskip

For (ii), note that $\chi(H)=|H|\langle{\chi_H,{\bf{1}}_H}\rangle$, where $\chi_H$ is the restriction of the character $\chi$ to $H$ and ${\bf{1}}_H$ is the restriction of the trivial character of $G$ to $H$. If $H\subseteq\ker(\chi)$, then $\chi(h)=\chi(1)$ for all $h\in H$. Thus $\chi(H)=|H|\chi(1)$. If $H\nsubseteq\ker(\chi)$, then by [\cite{isaacs1994character}, Corollary $6.7$], $\langle{\chi_H,{\bf{1}}_H}\rangle=0$ and the result follows.\qed

\section{Integrality of \texorpdfstring{$\Gamma_H(G)$}{Gamma	extunderscore H(G)}}\label{section-3}

In this section, we prove Theorem \ref{H-normal-graph-integral}, which provides a sufficient condition for the integrality of the graph $\Gamma_H(G)$. We first recall the notions of atoms and Eulerian subsets of a group $G$ and borrow an important result from \cite{guo2019integral}.

For $g\in G$, the \textit{atom} \cite{alperin2012integral} of $g$, denoted by $\operatorname{Atom}(g)$ is defined as $\lbrace{x\in G\mid \langle{x}\rangle=\langle{g}\rangle}\rbrace$, i.e., 
$$\operatorname{Atom}(g)=\lbrace{g^k\mid \operatorname{gcd}(k,o(g))=1}\rbrace.$$
A set $S\subseteq G$ is said to be \textit{Eulerian} \cite{guo2019integral} if for any $s\in S$, $\operatorname{Atom}(s)\subseteq S$. In other words, a set $S\subseteq G$ is Eulerian if $S$ is a (disjoint) union of some atoms in $G$.

\begin{proposition}\cite[Corollary 3]{guo2019integral}\label{Cay-int-suff-cond}
Let $G$ be a finite group and $H$ be a subgroup of $G$. Let $R$ be a normal, Eulerian subset of $G$. Then the graph $\operatorname{Cay}(G,S)$ is integral where $S\coloneqq R\setminus(R\cap H)$.
\end{proposition}

The following result is a routine exercise in group theory and therefore stated without a proof.

\begin{lemma}\label{S-normal-complement}
Let $G$ be a group and let $S\subseteq G$. Then $S$ is normal in $G$ if and only if $G\setminus S$ is normal in $G$.
\end{lemma}

For the vertices of Type 1,2 and 3 in $G\times G$, we define the following sets,
$$S_{H}^{(1)}\coloneqq\lbrace{(g,1)\mid g\in G\setminus H}\rbrace=(G\setminus H)\times\lbrace{1}\rbrace,$$
$$S_{H}^{(2)}\coloneqq\lbrace{(1,g)\mid g\in G\setminus H}\rbrace=\lbrace{1}\rbrace\times(G\setminus H),$$
$$S_{H}^{(3)}\coloneqq\lbrace{(g,g)\mid g\in G\setminus H}\rbrace=\Delta_{(G\setminus H)}=\Delta_G\setminus\Delta_H,$$
where $\Delta_X=\lbrace{(x,x)\mid x\in X}\rbrace$ denotes the diagonal of the set $X$. Define the graphs $\Gamma_{H}^{(k)}(G)\coloneqq\operatorname{Cay}\left(G\times G,S_{H}^{(k)}\right)$ for $k=1,2,3$.

\begin{lemma}\label{iso-1-2-3}
For $H\le G$, $\Gamma_{H}^{(1)}(G)\cong\Gamma_{H}^{(2)}(G)\cong\Gamma_{H}^{(3)}(G)$.
\end{lemma}
\pf Define the map $\alpha: G \times G \to G \times G$ as
$$(x,y)\mapsto(y^{-1},y^{-1}x)\quad \text{for all} \quad x,y \in G.$$
One can check that $\alpha$ is a graph isomorphism between $\Gamma_{H}^{(1)}(G)$ and $\Gamma_{H}^{(2)}(G)$ and also between $\Gamma_{H}^{(2)}(G)$ and $\Gamma_{H}^{(3)}(G)$. Hence  we have the result.\qed

\begin{lemma}\label{Gamma-1-2-3-integral}
For $H\le G$, the graphs $\Gamma_{H}^{(k)}(G)$ are integral for $k=1,2,3$.
\end{lemma}
\pf By Lemma \ref{iso-1-2-3}, it only suffices to show that $\Gamma_{H}^{(1)}(G)$ is integral. Consider the set $G\times\lbrace{1}\rbrace\subseteq G\times G$. Then it is clear that $G\times\lbrace{1}\rbrace$ is normal and Eulerian in $G\times G$. Consider the subgroup $H\times\lbrace{1}\rbrace$ of $G\times G$. Then a simple calculation yields that, $$(G\times\lbrace{1}\rbrace)\setminus((G\times\lbrace{1}\rbrace)\cap(H\times\lbrace{1}\rbrace))=S_{H}^{(1)}.$$
Hence $\Gamma_{H}^{(1)}(G)$ is integral by Proposition \ref{Cay-int-suff-cond}.\qed

\medskip

Now we are in a position to prove our main result of this section.

\begin{theorem}\label{H-normal-graph-integral}
Let $G$ be a group and $H\lhd G$. Then $\Gamma_H(G)$ is integral.
\end{theorem}
\pf For convenience, we do not distinguish between the representations of the group and their corresponding representation matrices throughout the proof. Let $\rho_1, \ldots, \rho_t$ be all inequivalent irreducible representations of $G$ with degrees $d_1, \ldots, d_t$ respectively. Then the set of inequivalent irreducible representations of the group $G\times G$ is given by 
$$\operatorname{IRR}(G\times G)=\lbrace{\Omega_{ij}\mid 1\le i,j\le t}\rbrace,$$
where $\Omega_{ij}\coloneqq\rho_i\otimes\rho_j$. Let $\mathcal{L}$ be the regular representation of $G\times G$. From Lemma \ref{lem:cayley1}, we have
\begin{equation}\label{regular-rep-equiv-tensor-product-rep}
\mathcal{L}\sim\bigoplus_{1\le i,j\le t}\Omega_{ij}^{\oplus d_id_j}.
\end{equation}
By Lemma \ref{lemma-babai}, the adjacency matrix of $\Gamma_H(G)$ becomes
$$\mathcal{A}(\Gamma_H(G))=\sum_{s\in S_H}\mathcal{L}(s)=\sum_{g\in G\setminus H}\mathcal{L}(g,1)+\sum_{g\in G\setminus H}\mathcal{L}(1,g)+\sum_{g\in G\setminus H}\mathcal{L}(g,g).$$
Thus, by \eqref{eqn:unique}, we obtain an orthogonal matrix $\mathcal{P}$ such that
{\small
\begin{equation}\label{eqn-2}
\mathcal{P}(\mathcal{A}(\Gamma_H(G)))\mathcal{P}^{-1}=\mathcal{P}\left(\sum_{s\in S_H}\mathcal{L}(s)\right)\mathcal{P}^{-1}=\bigoplus_{1\le i,j\le t}\left(\sum_{s\in S_H}\Omega_{ij}(s)\right)^{\oplus d_id_j}=\bigoplus_{1\le i,j\le t}\left(\Omega_{ij}(S_H)\right)^{\oplus d_id_j}.
\end{equation}}
It suffices to show that all the eigenvalues of $\Omega_{ij}(S_H)$ are integers. Let $\chi_i$ be the irreducible character of $G$ corresponding to the irreducible representation $\rho_i$. Then 
$$\operatorname{Irr}(G\times G)=\lbrace{\psi_{ij}\mid 1\le i,j\le t}\rbrace,$$
where $\psi_{ij}\coloneqq\chi_i\chi_j$. Now we can write $\Omega_{ij}(S_H)$ as,
\begin{equation}\label{Omega-ij-S_H}
\Omega_{ij}(S_H)=\Omega_{ij}\left(S_{H}^{(1)}\right)+\Omega_{ij}\left(S_{H}^{(2)}\right)+\Omega_{ij}\left(S_{H}^{(3)}\right).
\end{equation}
Since $H$ is normal in $G$, $G\setminus H$ is normal in $G$ (by Lemma \ref{S-normal-complement}) and hence both $S_{H}^{(1)}$ and $S_{H}^{(2)}$ are normal in $G\times G$. By Lemma \ref{chi-afford-rho-lemma}, we have
{\footnotesize\begin{equation}\label{Omega-S_H(1)-eqn}
\Omega_{ij}\left(S_H^{(1)}\right)=\frac{\psi_{ij}\left(S_{H}^{(1)}\right)}{\psi_{ij}(1,1)}I_{\psi_{ij}(1,1)}=\frac{(\chi_i\chi_j)\left(S_{H}^{(1)}\right)}{(\chi_i\chi_j)(1,1)}I_{(\chi_i\chi_j)(1,1)}=\frac{\chi_i(G\setminus H)\chi_j(1)}{\chi_i(1)\chi_j(1)}I_{\chi_i(1)\chi_j(1)}=\frac{\chi_i(G\setminus H)}{d_i}I_{d_id_j}.
\end{equation}}
{\footnotesize\begin{equation}\label{Omega-S_H(2)-eqn}
\Omega_{ij}\left(S_H^{(2)}\right)=\frac{\psi_{ij}\left(S_{H}^{(2)}\right)}{\psi_{ij}(1,1)}I_{\psi_{ij}(1,1)}=\frac{(\chi_i\chi_j)\left(S_{H}^{(2)}\right)}{(\chi_i\chi_j)(1,1)}I_{(\chi_i\chi_j)(1,1)}=\frac{\chi_i(1)\chi_j(G\setminus H)}{\chi_i(1)\chi_j(1)}I_{\chi_i(1)\chi_j(1)}=\frac{\chi_j(G\setminus H)}{d_j}I_{d_id_j}.
\end{equation}}
Thus from \eqref{Omega-ij-S_H}, we have
\begin{equation}\label{Omega-S_H}
\Omega_{ij}(S_H)=\frac{\chi_i(G\setminus H)}{d_i}I_{d_id_j}+\frac{\chi_j(G\setminus H)}{d_j}I_{d_id_j}+\Omega_{ij}\left(S_{H}^{(3)}\right).
\end{equation}
Note that by Lemma \ref{lemma-babai}, the adjacency matrix of $\Gamma_{H}^{(k)}(G)$ for $k=1,2,3$ becomes
$$\mathcal{A}\left(\Gamma_{H}^{(1)}(G)\right)=\sum_{s\in S_{H}^{(1)}}\mathcal{L}(s)=\sum_{g\in G\setminus H}\mathcal{L}(g,1),$$
$$\mathcal{A}\left(\Gamma_{H}^{(2)}(G)\right)=\sum_{s\in S_{H}^{(2)}}\mathcal{L}(s)=\sum_{g\in G\setminus H}\mathcal{L}(1,g),$$
$$\mathcal{A}\left(\Gamma_{H}^{(3)}(G)\right)=\sum_{s\in S_{H}^{(3)}}\mathcal{L}(s)=\sum_{g\in G\setminus H}\mathcal{L}(g,g).$$
Hence for each $k=1,2,3$, we have the same orthogonal matrix $\mathcal{P}$ satisfying
$$\mathcal{P}(\mathcal{A}(\Gamma_{H}^{(k)}(G)))\mathcal{P}^{-1}=\bigoplus_{1\le i,j\le t}\left(\sum_{s\in S_{H}^{(k)}}\Omega_{ij}(s)\right)^{\oplus d_id_j}=\bigoplus_{1\le i,j\le t}\left(\Omega_{ij}\left(S_H^{(k)}\right)\right)^{\oplus d_id_j}.$$
Since $\Gamma_{H}^{(k)}(G)$, $k=1,2,3$ are integral by Lemma \ref{Gamma-1-2-3-integral}, we have that all the eigenvalues of $\Omega_{ij}\left(S_{H}^{(k)}\right)$ are integers. Since $\Omega_{ij}\left(S_H^{(1)}\right)$ and $\Omega_{ij}\left(S_H^{(2)}\right)$ are scalar matrices, from \eqref{Omega-S_H(1)-eqn} and \eqref{Omega-S_H(2)-eqn}, we have $\frac{\chi_i(G\setminus H)}{d_i}$ and $\frac{\chi_j(G\setminus H)}{d_j}$ are integers for all $1\le i,j\le t$. Consequently from \eqref{Omega-S_H} all the eigenvalues of $\Omega_{ij}(S_H)$ are integers and hence $\Gamma_H(G)$ is integral by \eqref{eqn-2}.\qed

\section{A few eigenvalues of \texorpdfstring{$\Gamma_H(G)$}{Gamma\_H(G)}}\label{section-4}

In this section, we obtain some general results regarding the eigenvalues of $\Gamma_H(G)$. First we obtain a few (repeated) eigenvalues of $\Gamma_H(G)$ and then obtain a sufficient condition for $\Gamma_H(G)$ to be singular. For this purpose, we recall two folklore results on groups and Cayley graphs.

\begin{lemma}\label{subgroup-complement-gen-G}
Let $G$ be a group and $H\le G$. Then $G\setminus H$ generates $G$.
\end{lemma}

\begin{lemma}\label{Cay-cartesian-prod-lemma}
Let $G_1,G_2$ be two groups and let $S_1\subseteq G_1\setminus\lbrace{1_{G_1}}\rbrace,S_2\subseteq G_2\setminus\lbrace{1_{G_2}}\rbrace$ be such that $S_1^{-1}=S_1$ and $S_2^{-1}=S_2$. Then
$$\operatorname{Cay}(G_1\times G_2,(S_1\times\lbrace{1_{G_2}}\rbrace\cup\lbrace{1_{G_1}}\rbrace\times S_2))\cong\operatorname{Cay}(G_1,S_1)\square\operatorname{Cay}(G_2,S_2)$$
where $\square$ denotes the cartesian product of graphs.
\end{lemma}

\begin{proposition}\label{simple-eigenval}
$3(|G|-|H|)$ is a simple eigenvalue of $\Gamma_H(G)$.
\end{proposition}
\pf Let $T\coloneqq S_H^{(1)}\cup S_H^{(2)}$. Then by Lemma \ref{Cay-cartesian-prod-lemma}, 
$$\operatorname{Cay}(G\times G,T)\cong\operatorname{Cay}(G,G\setminus H)\square\operatorname{Cay}(G,G\setminus H).$$
By Lemma \ref{subgroup-complement-gen-G}, $\operatorname{Cay}(G,G\setminus H)$ is connected and hence $\operatorname{Cay}(G\times G,T)$ is connected, which implies the graph $\Gamma_H(G)$ is connected. Since $\Gamma_H(G)$ is $|S_H|$-regular, by Perron-Frobenius theorem, $3(|G|-|H|)$ is a simple eigenvalue of $\Gamma_H(G)$.\qed

\begin{lemma}\label{eigenvalue-projection-operator-irr-lemma}
Let $G$ be a group and $H\le G$. Let $(\rho,V)\in\operatorname{IRR}(G)$. Then
\begin{enumerate}[label=(\roman*)]
\item The eigenvalue of $\rho(G)\coloneqq\sum_{g\in G}\rho(g)$ is $|G|$ if $\rho$ is the trivial representation and $0$ otherwise. More specifically, $\rho(G)=0$, if $\rho$ is a non-trivial representation.
\item The eigenvalues of $\rho(H)\coloneqq\sum_{h\in H}\rho(h)$ are
$$\left[\underbrace{|H|,|H|,\ldots,|H|}_{\dim\left(V^H\right)~times},\underbrace{0,0,\ldots,0}_{\dim V-\dim\left(V^H\right)~times}\right].$$
\end{enumerate}
\end{lemma}
\pf (i) directly follows from Lemma \ref{chi-afford-rho-lemma} and Lemma \ref{G-H-char-lemma} (i). For (ii), note that 
$$P_H\coloneqq\frac{1}{|H|}\sum_{h\in H}\rho(h)$$
is a projection from $V$ onto its $H$-fixed subspace $V^H$. Since a projection operator is diagonalizable, 
$$\#\text{ Non-zero eigenvalues of } P_H=\operatorname{Rank}(P_H)=\dim\left(V^H\right).$$
Also, the only eigenvalues of a projection operator are $0$ and $1$. Hence we obtain the eigenvalues of $P_H$ as
$$\left[\underbrace{1,1,\ldots,1}_{\dim\left(V^H\right)~\text{times}},\underbrace{0,0,\ldots,0}_{\dim V-\dim\left(V^H\right)~\text{times}}\right].$$
This completes the proof.\qed

\begin{theorem}\label{other-two-eigenval}
Let $G$ be a group of order $n$ and $H\le G$. Then 
$$|G|-|H|,|G|-3|H|\in \operatorname{Spec}(\Gamma_H(G)).$$ Moreover, if $[G:H]=\ell$, then $m(|G|-3|H|)\ge 2(\ell-1)$ and $m(|G|-|H|)\ge 2(n-\ell)$.
\end{theorem}
\pf Let 
$$\operatorname{IRR}(G)=\lbrace{(\rho_i,V_i)\mid 1\le i\le t}\rbrace,$$ 
and 
$$\operatorname{Irr}(G)=\lbrace{\chi_i\mid 1\le i\le t}\rbrace,$$
where $\rho_1$ is the trivial representation and $\chi_1$ is the trivial character of $G$. Let $d_i=\deg(\chi_i)=\deg(\rho_i)$. In \eqref{Omega-ij-S_H}, we take $i=1$ and $j\ge 2$. Then we have
\begin{align*}
\Omega_{1j}(S_H)&=\Omega_{1j}\left(S_{H}^{(1)}\right)+\Omega_{1j}\left(S_{H}^{(2)}\right)+\Omega_{1j}\left(S_{H}^{(3)}\right)\\
&=(|G|-|H|)I_{d_j}+2\left(\sum_{g\in G\setminus H}\rho_j(g)\right).
\end{align*}
For $2\le j\le t$, $\rho_j(G\setminus H)=-\rho_j(H)$, by Lemma \ref{eigenvalue-projection-operator-irr-lemma} (i). Thus from Lemma \ref{eigenvalue-projection-operator-irr-lemma} (ii), for $j\ge 2$, we have the eigenvalues of $\rho_j(G\setminus H)$ are 
$$\left[\underbrace{-|H|,-|H|,\ldots,-|H|}_{\dim\left(V_j^H\right)~\text{times}},\underbrace{0,0,\ldots,0}_{\dim V_j-\dim\left(V_j^H\right)~\text{times}}\right].$$
Thus for $2\le j\le t$, the eigenvalues of $\Omega_{1j}(S_H)$ are 
$$\left[\underbrace{|G|-3|H|,|G|-3|H|,\ldots,|G|-3|H|}_{\dim\left(V_j^H\right)~\text{times}},\underbrace{|G|-|H|,|G|-|H|,\ldots,|G|-|H|}_{\dim V_j-\dim\left(V_j^H\right)~\text{times}}\right].$$
This proves the first part.

To compute the multiplicity, note that for $2\le i\le t$, a similar approach yields the eigenvalues of $\Omega_{i1}(S_H)$ as 
$$\left[\underbrace{|G|-3|H|,|G|-3|H|,\ldots,|G|-3|H|}_{\dim\left(V_i^H\right)~\text{times}},\underbrace{|G|-|H|,|G|-|H|,\ldots,|G|-|H|}_{\dim V_i-\dim\left(V_i^H\right)~\text{times}}\right].$$
From \eqref{regular-rep-equiv-tensor-product-rep}, for each $2\le j\le t$, $|G|-3|H|$ occurs $d_1d_j\dim\left(V_j^H\right)$ times in $\Gamma_H(G)$ and for each $2\le i\le t$, $|G|-3|H|$ occurs $d_id_1\dim\left(V_i^H\right)$ times in $\Gamma_H(G)$. Similarly, for each $2\le j\le t$, $|G|-|H|$ occurs $d_1d_j\left(\dim V_j-\dim\left(V_j^H\right)\right)$ times in $\Gamma_H(G)$ and for each $2\le i\le t$, $|G|-|H|$ occurs $d_id_1\left(\dim V_i-\dim\left(V_i^H\right)\right)$ times in $\Gamma_H(G)$. Thus
\begin{align*}
m(|G|-3|H|)&\ge\sum_{2\le j\le t}d_1d_j\dim\left(V_j^H\right)+\sum_{2\le i\le t}d_id_1\dim\left(V_i^H\right)\\
&=2\left(\sum_{i=2}^td_i\dim\left(V_i^H\right)\right),
\end{align*}
and
\begin{align*}
m(|G|-|H|)&\ge\sum_{2\le j\le t}d_1d_j\left(\dim V_j-\dim\left(V_j^H\right)\right)+\sum_{2\le i\le t}d_id_1\left(\dim V_i-\dim\left(V_i^H\right)\right)\\
&=2\left(\sum_{i=2}^td_i\left(\dim V_i-\dim\left(V_i^H\right)\right)\right)\\
&=2\left(\sum_{i=2}^td_i^2-\sum_{i=2}^td_i\dim\left(V_i^H\right)\right).
\end{align*}
Now
\begin{align*}
\sum_{i=2}^td_i\dim (V_i^H)&=\sum_{i=2}^td_i\left(\frac{1}{|H|}\sum_{h\in H}\chi_i(h)\right)\\
&=\frac{1}{|H|}\sum_{h\in H}\left(\sum_{i=2}^td_i\chi_i(h)\right)\\
&=\frac{1}{|H|}\sum_{h\in H}\left(\sum_{i=1}^td_i\chi_i(h)-d_1\chi_1(h)\right)\\
&=\frac{1}{|H|}\sum_{h\in H}(\chi_L(h)-1),
\end{align*}
where $\chi_L$ is the character of the regular representation. Hence,
\begin{align*}
\sum_{i=2}^td_i\dim\left(V_i^H\right)&=\frac{1}{|H|}\left[(\chi_L(1)-1)+\sum_{h\in H,h\ne 1}(\chi_L(h)-1)\right]\\
&=\frac{1}{|H|}[(n-1)-(|H|-1)].
\end{align*}
Putting $|H|=n/\ell$, we have $$\sum_{i=2}^td_i\dim\left(V_i^H\right)=\ell-1.$$
Thus 
\begin{align*}
m(|G|-3|H|)&\ge 2\left(\sum_{i=2}^td_i\dim\left(V_i^H\right)\right)\\
&=2(\ell-1),
\end{align*}
and 
\begin{align*}
m(|G|-|H|)&\ge 2\left(\sum_{i=2}^td_i^2-\sum_{i=2}^td_i\dim\left(V_i^H\right)\right)\\
&=2(n-\ell). 
\end{align*}
\qed

\begin{theorem}
 Let $H\le G$ with $|H|>2$. Then $\Gamma_H(G)$ is singular.
\end{theorem}

\pf
We construct an eigenvector $v$ such that
$\mathcal{A}(\Gamma_H(G))v=\bf{0}$ holds.
As $|H|>2$, our strategy is to take two distinct elements
$h,h'\in H\setminus\lbrace 1\rbrace$ and relabel the vertices
accordingly. We divide the proof into two cases and show that
in both cases, the vector
$$
v=(1,1,1,-1,-1,-1,
\underbrace{0,0,\ldots,0}_{(n^2-6)~\text{times}})^T
$$
serves our purpose.

\medskip

\textbf{Case 1:}
Suppose $H$ contains an element $h$ of order $\ge 3$,
i.e., $h\neq h^{-1}$. Choose $h'=h^{-1}$.
We relabel the vertices $a_1,\ldots,a_{n^2}\in G\times G$
as
$a_1=(h,1)$, $a_2=(1,h)$, $a_3=(h',h')$,
$a_4=(h',1)$, $a_5=(1,h')$, $a_6=(h,h)$
in this case.

Let $(g_i,g_j)\in G\times G$.
If $(g_i,g_j)$ is adjacent to none of
$(h,1),(1,h),(h',h'),(h',1),(1,h'),(h,h)$,
then the first six entries of the row in
$\mathcal{A}(\Gamma_H(G))$ corresponding to $(g_i,g_j)$,
say $R_{(g_i,g_j)}$, are $0,0,0,0,0,0$, i.e.,
$R_{(g_i,g_j)}=(0,0,0,0,0,0,*,\ldots,*)$.
Hence $R_{(g_i,g_j)}v=0$.

\medskip

\textbf{Subcase 1:}
Suppose that $(g_i,g_j)\sim(h,1)$, i.e.,
$(g_ih^{-1},g_j)\in S_H$.
Then the following cases can occur:

\begin{itemize}

\item
If $(g_ih^{-1},g_j)\in S_H^{(1)}$, then
$g_ih^{-1}\in G\setminus H$ and $g_j=1$.
Hence $(g_i,g_j)\sim(h',1)$ and
$(g_i,g_j)\nsim(1,h),(h',h'),(h,h),(1,h')$.
Therefore
$R_{(g_i,g_j)}=(1,0,0,1,0,0,*,\ldots,*)$,
i.e., $R_{(g_i,g_j)}v=0$.
The matrix form corresponding to this case is shown below:

$$
{\footnotesize
\begin{pNiceMatrix}[first-row=\scriptstyle,first-col=\scriptstyle]
& (h,1) & (1,h) & (h',h') & (h',1) & (1,h') & (h,h)\\
(h,1) & 0 & 0 & 0 & 0 & 0 & 0 & * & \cdots & *\\
(1,h) & 0 & 0 & 0 & 0 & 0 & 0 & * & \cdots & *\\
(h',h') & 0 & 0 & 0 & 0 & 0 & 0 & * & \cdots & *\\
(h',1) & 0 & 0 & 0 & 0 & 0 & 0 & * & \cdots & *\\
(1,h') & 0 & 0 & 0 & 0 & 0 & 0 & * & \cdots & *\\
(h,h) & 0 & 0 & 0 & 0 & 0 & 0 & * & \cdots & *\\
 & * & * & * & * & * & * & * & \cdots & *\\
 & \vdots & \vdots & \vdots & \vdots & \vdots
 & \vdots & \vdots & \ddots & \vdots\\
 & * & * & * & * & * & * & * & \cdots & *\\
(g_i,g_j) & 1 & 0 & 0 & 1 & 0 & 0 & * & \cdots & *\\
 & * & * & * & * & * & * & * & \cdots & *\\
 & \vdots & \vdots & \vdots & \vdots & \vdots
 & \vdots & \vdots & \ddots & \vdots\\
 & * & * & * & * & * & * & * & \cdots & *
\end{pNiceMatrix}}
{\begin{pmatrix}
1\\1\\1\\-1\\-1\\-1\\0\\\vdots\\0
\end{pmatrix}=\bf{0}}.
$$

\item
If $(g_ih^{-1},g_j)\in S_H^{(2)}$, then
$g_ih^{-1}=1$ and $g_j\in G\setminus H$.
Hence $(g_i,g_j)\sim(h,h)$ and
$(g_i,g_j)\nsim(1,h),(h',h'),(h',1),(1,h')$.
Therefore
$R_{(g_i,g_j)}=(1,0,0,0,0,1,*,\ldots,*)$,
i.e., $R_{(g_i,g_j)}v=0$.

\item
If $(g_ih^{-1},g_j)\in S_H^{(3)}$, then
$g_ih^{-1}=g_j\in G\setminus H$.
Hence $(g_i,g_j)\sim(1,h')$ and
$(g_i,g_j)\nsim(1,h),(h',h'),(h',1),(h,h)$.
Therefore
$R_{(g_i,g_j)}=(1,0,0,0,1,0,*,\ldots,*)$,
i.e., $R_{(g_i,g_j)}v=0$.

\end{itemize}

\textbf{Subcase 2:}
Suppose that $(g_i,g_j)\sim(1,h)$, i.e.,
$(g_i,g_jh^{-1})\in S_H$.
Then the following cases can occur:

\begin{itemize}

\item
If $(g_i,g_jh^{-1})\in S_H^{(1)}$, then
$g_i\in G\setminus H$ and $g_jh^{-1}=1$.
Hence $(g_i,g_j)\sim(h,h)$ and
$(g_i,g_j)\nsim(h,1),(h',h'),(h',1),(1,h')$.
Therefore
$R_{(g_i,g_j)}=(0,1,0,0,0,1,*,\ldots,*)$,
i.e., $R_{(g_i,g_j)}v=0$.

\item
If $(g_i,g_jh^{-1})\in S_H^{(2)}$, then
$g_i=1$ and $g_jh^{-1}\in G\setminus H$.
Hence $(g_i,g_j)\sim(1,h')$ and
$(g_i,g_j)\nsim(h,1),(h',h'),(h,h),(h',1)$.
Therefore
$R_{(g_i,g_j)}=(0,1,0,0,1,0,*,\ldots,*)$,
i.e., $R_{(g_i,g_j)}v=0$.

\item
If $(g_i,g_jh^{-1})\in S_H^{(3)}$, then
$g_i=g_jh^{-1}\in G\setminus H$.
Hence $(g_i,g_j)\sim(h',1)$ and
$(g_i,g_j)\nsim(h,1),(h',h'),(1,h'),(h,h)$.
Therefore
$R_{(g_i,g_j)}=(0,1,0,1,0,0,*,\ldots,*)$,
i.e., $R_{(g_i,g_j)}v=0$.

\end{itemize}

\textbf{Subcase 3:}
Suppose that $(g_i,g_j)\sim(h,h)$, i.e.,
$(g_ih^{-1},g_jh^{-1})\in S_H$.
Then the following cases can occur:

\begin{itemize}

\item
If $(g_ih^{-1},g_jh^{-1})\in S_H^{(1)}$, then
$g_ih^{-1}\in G\setminus H$ and $g_jh^{-1}=1$.
Hence $(g_i,g_j)\sim(1,h)$ and
$(g_i,g_j)\nsim(h,1),(h',h'),(h',1),(1,h')$.
Therefore
$R_{(g_i,g_j)}=(0,1,0,0,0,1,*,\ldots,*)$,
i.e., $R_{(g_i,g_j)}v=0$.

\item
If $(g_ih^{-1},g_jh^{-1})\in S_H^{(2)}$, then
$g_ih^{-1}=1$ and
$g_jh^{-1}\in G\setminus H$.
Hence $(g_i,g_j)\sim(h,1)$ and
$(g_i,g_j)\nsim(1,h),(h',h'),(h',1),(1,h')$.
Therefore
$R_{(g_i,g_j)}=(1,0,0,0,0,1,*,\ldots,*)$,
i.e., $R_{(g_i,g_j)}v=0$.

\item If $(g_ih^{-1},g_jh^{-1})\in S_H^{(3)}$, then $g_ih^{-1}=g_jh^{-1}\in G\setminus H$.
Hence $(g_i,g_j)\sim(h',h')$ and
$(g_i,g_j)\nsim(h,1),(1,h),(h',1),(1,h')$. Therefore $R_{(g_i,g_j)}=(0,0,1,0,0,1,*,\ldots,*)$, i.e., $R_{(g_i,g_j)}v=0$.

\end{itemize}

The remaining cases, where $(g_i,g_j)$ is adjacent to one of $(h',h')$, $(h',1)$, or $(1,h')$, are summarized in the following table. Here, $k$ indicates that $(g_i,g_j)u^{-1}\in S_H^{(k)}$. The last column records the first six entries of $R_{(g_i,g_j)}$.

$$
\begin{array}{|c|c|c|c|}
\hline
u & k & \textbf{Other adjacent vertex}
& \textbf{First six entries of }R_{(g_i,g_j)}
\\ \hline
(h',h') & 1 & (1,h')  & (0,0,1,0,1,0)\\
        & 2 & (h',1)  & (0,0,1,1,0,0)\\
        & 3 & (h,h)   & (0,0,1,0,0,1)\\ \hline
(h',1)  & 1 & (h,1)   & (1,0,0,1,0,0)\\
        & 2 & (h',h') & (0,0,1,1,0,0)\\
        & 3 & (1,h)   & (0,1,0,1,0,0)\\ \hline
(1,h')  & 1 & (h',h') & (0,0,1,0,1,0)\\
        & 2 & (1,h)   & (0,1,0,0,1,0)\\
        & 3 & (h,1)   & (1,0,0,0,1,0)\\
\hline
\end{array}
$$

In each row, exactly one of the first three entries and exactly one of the next three entries are equal to $1$. Thus $R_{(g_i,g_j)}v=0$ in every remaining case. Hence every row of $\mathcal{A}(\Gamma_H(G))$ is orthogonal to $v$, and therefore $\mathcal{A}(\Gamma_H(G))v=\mathbf{0}$.

\medskip

\textbf{Case 2.}
Suppose every non-identity element of $H$ has order $2$.
Therefore $H$ is abelian.
Choose distinct elements
$h,h'\in H\setminus\{1\}$.
We relabel the vertices
$a_1,\ldots,a_{n^2}\in G\times G$
as
$a_1=(h,1)$, $a_2=(1,h')$, $a_3=(h',h)$,
$a_4=(h',1)$, $a_5=(1,h)$, $a_6=(h,h')$
in this case.

Let $(g_i,g_j)\in G\times G$.
If $(g_i,g_j)$ is adjacent to none of
$(h,1),(1,h'),(h',h),(h',1),(1,h),(h,h')$,
then the first six entries of the row in
$\mathcal{A}(\Gamma_H(G))$ corresponding to
$(g_i,g_j)$, say $R_{(g_i,g_j)}$,
are $0,0,0,0,0,0$, i.e.,
$R_{(g_i,g_j)}=(0,0,0,0,0,0,*,\ldots,*)$.
Hence $R_{(g_i,g_j)}v=0$.

\medskip

\textbf{Subcase 1:}
Suppose that $(g_i,g_j)\sim(h,1)$, i.e.,
$(g_ih^{-1},g_j)\in S_H$.
Then the following cases can occur:

\begin{itemize}

\item
If $(g_ih^{-1},g_j)\in S_H^{(1)}$, then
$g_ih^{-1}\in G\setminus H$ and $g_j=1$.
Hence $(g_i,g_j)\sim(h',1)$ and
$(g_i,g_j)\nsim(1,h'),(h',h),(1,h),(h,h')$.
Therefore
$R_{(g_i,g_j)}=(1,0,0,1,0,0,*,\ldots,*)$,
i.e., $R_{(g_i,g_j)}v=0$.
The matrix form corresponding to this case is shown below:

$$
{\footnotesize
\begin{pNiceMatrix}[first-row=\scriptstyle,first-col=\scriptstyle]
& (h,1) & (1,h') & (h',h) & (h',1) & (1,h) & (h,h')\\
(h,1) & 0 & 0 & 0 & 0 & 0 & 0 & * & \cdots & *\\
(1,h') & 0 & 0 & 0 & 0 & 0 & 0 & * & \cdots & *\\
(h',h) & 0 & 0 & 0 & 0 & 0 & 0 & * & \cdots & *\\
(h',1) & 0 & 0 & 0 & 0 & 0 & 0 & * & \cdots & *\\
(1,h) & 0 & 0 & 0 & 0 & 0 & 0 & * & \cdots & *\\
(h,h') & 0 & 0 & 0 & 0 & 0 & 0 & * & \cdots & *\\
 & * & * & * & * & * & * & * & \cdots & *\\
 & \vdots & \vdots & \vdots & \vdots & \vdots
 & \vdots & \vdots & \ddots & \vdots\\
 & * & * & * & * & * & * & * & \cdots & *\\
(g_i,g_j) & 1 & 0 & 0 & 1 & 0 & 0 & * & \cdots & *\\
 & * & * & * & * & * & * & * & \cdots & *\\
 & \vdots & \vdots & \vdots & \vdots & \vdots
 & \vdots & \vdots & \ddots & \vdots\\
 & * & * & * & * & * & * & * & \cdots & *
\end{pNiceMatrix}}
{\begin{pmatrix}
1\\1\\1\\-1\\-1\\-1\\0\\\vdots\\0
\end{pmatrix}=\bf{0}}.
$$

\item
If $(g_ih^{-1},g_j)\in S_H^{(2)}$, then
$g_ih^{-1}=1$ and $g_j\in G\setminus H$.
Hence $(g_i,g_j)\sim(h,h')$ and
$(g_i,g_j)\nsim(1,h'),(h',h),(h',1),(1,h)$.
Therefore
$R_{(g_i,g_j)}=(1,0,0,0,0,1,*,\ldots,*)$,
i.e., $R_{(g_i,g_j)}v=0$.

\item
If $(g_ih^{-1},g_j)\in S_H^{(3)}$, then
$g_ih^{-1}=g_j\in G\setminus H$.
Hence $(g_i,g_j)\sim(1,h)$ and
$(g_i,g_j)\nsim(1,h'),(h',h),(h',1),(h,h')$.
Therefore
$R_{(g_i,g_j)}=(1,0,0,0,1,0,*,\ldots,*)$,
i.e., $R_{(g_i,g_j)}v=0$.

\end{itemize}

\textbf{Subcase 2:}
Suppose that $(g_i,g_j)\sim(1,h)$, i.e.,
$(g_i,g_jh^{-1})\in S_H$.
Then the following cases can occur:

\begin{itemize}

\item
If $(g_i,g_jh^{-1})\in S_H^{(1)}$, then
$g_i\in G\setminus H$ and $g_jh^{-1}=1$.
Hence $(g_i,g_j)\sim(h',h)$ and
$(g_i,g_j)\nsim(h,1),(1,h'),(h',1),(h,h')$.
Therefore
$R_{(g_i,g_j)}=(0,0,1,0,1,0,*,\ldots,*)$,
i.e., $R_{(g_i,g_j)}v=0$.

\item
If $(g_i,g_jh^{-1})\in S_H^{(2)}$, then
$g_i=1$ and $g_jh^{-1}\in G\setminus H$.
Hence $(g_i,g_j)\sim(1,h')$ and
$(g_i,g_j)\nsim(h,1),(h',h),(h',1),(h,h')$.
Therefore
$R_{(g_i,g_j)}=(0,1,0,0,1,0,*,\ldots,*)$,
i.e., $R_{(g_i,g_j)}v=0$.

\item
If $(g_i,g_jh^{-1})\in S_H^{(3)}$, then
$g_i=g_jh^{-1}\in G\setminus H$.
Hence $(g_i,g_j)\sim(h,1)$ and
$(g_i,g_j)\nsim(1,h'),(h',h),(h',1),(h,h')$.
Therefore
$R_{(g_i,g_j)}=(1,0,0,0,1,0,*,\ldots,*)$,
i.e., $R_{(g_i,g_j)}v=0$.

\end{itemize}

\textbf{Subcase 3:}
Suppose that $(g_i,g_j)\sim(h',h)$, i.e.,
$(g_i(h')^{-1},g_jh^{-1})\in S_H$.
Then the following cases can occur:

\begin{itemize}

\item
If $(g_i(h')^{-1},g_jh^{-1})\in S_H^{(1)}$,
then $g_i(h')^{-1}\in G\setminus H$
and $g_jh^{-1}=1$.
Hence $(g_i,g_j)\sim(1,h)$ and
$(g_i,g_j)\nsim(h,1),(1,h'),(h',1),(h,h')$.
Therefore
$R_{(g_i,g_j)}=(0,0,1,0,1,0,*,\ldots,*)$,
i.e., $R_{(g_i,g_j)}v=0$.

\item
If $(g_i(h')^{-1},g_jh^{-1})\in S_H^{(2)}$,
then $g_i(h')^{-1}=1$ and
$g_jh^{-1}\in G\setminus H$.
Hence $(g_i,g_j)\sim(h',1)$ and
$(g_i,g_j)\nsim(h,1),(1,h'),(1,h),(h,h')$.
Therefore
$R_{(g_i,g_j)}=(0,0,1,1,0,0,*,\ldots,*)$,
i.e., $R_{(g_i,g_j)}v=0$.

\item If $(g_i(h')^{-1},g_jh^{-1})\in S_H^{(3)}$,
then
$g_i(h')^{-1}=g_jh^{-1}\in G\setminus H$.
Hence $(g_i,g_j)\sim(h,h')$ as $hh'=h'h$
and it can be shown that
$(g_i,g_j)\nsim(h,1),(1,h'),(h',1),(1,h)$.
Therefore
$R_{(g_i,g_j)}=(0,0,1,0,0,1,*,\ldots,*)$,
i.e., $R_{(g_i,g_j)}v=0$.

\end{itemize}

The remaining cases, where $(g_i,g_j)$ is adjacent to one of $(h',1)$, $(1,h')$, or $(h,h')$, are summarized in the following table. Here, $k$ indicates that $(g_i,g_j)u^{-1}\in S_H^{(k)}$. The last column records the first six entries of $R_{(g_i,g_j)}$.

$$
\begin{array}{|c|c|c|c|}
\hline
u & k & \textbf{Other adjacent vertex}
& \textbf{First six entries of }R_{(g_i,g_j)}
\\ \hline
(h',1) & 1 & (h,1)   & (1,0,0,1,0,0)\\
       & 2 & (h',h)  & (0,0,1,1,0,0)\\
       & 3 & (1,h')  & (0,1,0,1,0,0)\\ \hline
(1,h') & 1 & (h,h')  & (0,1,0,0,0,1)\\
       & 2 & (1,h)   & (0,1,0,0,1,0)\\
       & 3 & (h',1)  & (0,1,0,1,0,0)\\ \hline
(h,h') & 1 & (1,h')  & (0,1,0,0,0,1)\\
       & 2 & (h,1)   & (1,0,0,0,0,1)\\
       & 3 & (h',h)  & (0,0,1,0,0,1)\\
\hline
\end{array}
$$

Again, each row has exactly one $1$ among the first three entries and exactly one $1$ among the next three entries. Therefore $R_{(g_i,g_j)}v=0$ in every remaining case. Hence every row of $\mathcal{A}(\Gamma_H(G))$ is orthogonal to $v$, and consequently
$\mathcal{A}(\Gamma_H(G))v=\mathbf{0}$. This completes the proof.
\qed

\section{Complete spectrum of \texorpdfstring{$\Gamma_H(G)$}{Gamma	extunderscore H(G)}, when \texorpdfstring{$H\lhd G$}{H normal G}}
\label{section-5}

In this section, we compute the complete spectrum of $\Gamma_H(G)$ for a group $G$ and its normal subgroup $H$, which essentially proves Theorem \ref{qsrg-isospectral}. A summary of the complete spectrum in various cases can be found in Table \ref{tab:summary-table-eigenvalues}.

\begin{table}[h]
\scriptsize
\centering
\renewcommand{\arraystretch}{1.3}
\begin{tabular}{|c|c|c|c|l|c|}
\hline
\textbf{$[G:H]$} & \textbf{$\kappa(\Gamma_H(G))$} & \textbf{Eigenvalues ($\lambda$)} & \textbf{Multiplicities ($m$)} & \textbf{\quad \quad \quad Source Blocks} & \textbf{Figure} \\ \hline
\multirow{4}{*}{\textbf{$2$}} & \multirow{4}{*}{$4$} & $3k$ & $1$ & \quad \quad \quad \quad \quad \quad $\{3k\}$ & \multirow{4}{*}{Fig. \ref{fig:Spec-[G:H]=2}} \\ 
 & & $k$ & $3n-6$ & \quad \quad $\mathcal{Q}_{i,j}$, and $\mathcal{U}$ inside $\mathcal{B}^{(2)}$ & \\ 
 & & $-k$ & $3n-3$ & $\{-k\}$, $\mathcal{R}_{i,j}$, and $\mathcal{U},\mathcal{V}$ inside $\mathcal{B}^{(2)}$ & \\ 
 & & $0$ & $n^2-6n+8$ & \quad \quad \quad $\mathcal{U}$ and $\mathcal{V}$ inside $\mathcal{B}^{(2)}$ & \\ \hline
\multirow{5}{*}{\textbf{$3$}} & \multirow{5}{*}{$5$} & $6k$ & $1$ & \quad \quad \quad \quad \quad \quad $\{6k\}$ & \multirow{5}{*}{Fig. \ref{fig:Spec-[G:H]=3}} \\ 
 & & $2k$ & $3n-9$ & \quad~ $\mathcal{Q}_{i,j}$, and $\mathcal{U}$ inside $\mathcal{B}^{(2)}$ & \\ 
 & & $-k$ & $6n-18$ & \quad~ $\mathcal{R}_{i,j}$, and $\mathcal{U},\mathcal{V}$ inside $\mathcal{B}^{(2)}$ & \\ 
 & & $-3k$ & $2$ & \quad \quad \quad \quad \quad~ $\{-3k\}$ & \\ 
 & & $0$ & $n^2-9n+24$ & \quad~ $\{0\}$, and $\mathcal{U},\mathcal{V}$ inside $\mathcal{B}^{(2)}$ & \\ \hline
\multirow{5}{*}{\makecell{\textbf{$\ell \ge 4$} \\ $(|H|=2)$}} & \multirow{5}{*}{$5$} & $6\ell-6$ & $1$ & \quad \quad \quad \quad~ $\{6(\ell-1)\}$ & \multirow{5}{*}{Fig. \ref{fig:Spec-[G:H]=l,k=2}} \\ 
 & & $2\ell-6$ & $3\ell-3$ & \quad \quad $\mathcal{P}_{i,j}$, and $\mathcal{S}$ inside $\mathcal{B}^{(1)}$ & \\ 
 & & $2\ell-2$ & $3\ell$ & \quad \quad $\mathcal{Q}_{i,j}$, and $\mathcal{U}$ inside $\mathcal{B}^{(2)}$ & \\ 
 & & $-6$ & $\ell^2-3\ell+2$ & \quad \quad~ $\mathcal{S}$ and $\mathcal{T}$ inside $\mathcal{B}^{(1)}$ & \\ 
 & & $-2$ & $3\ell^2-3\ell$ & \quad~ $\mathcal{R}_{i,j}$, and $\mathcal{U},\mathcal{V}$ inside $\mathcal{B}^{(2)}$ & \\ \hline
\multirow{6}{*}{\makecell{\textbf{$\ell \ge 4$} \\ $(|H|>2)$}} & \multirow{6}{*}{$6$} & $3(\ell-1)k$ & $1$ & \quad \quad \quad \quad~ $\{3(\ell-1)k\}$ & \multirow{6}{*}{Fig. \ref{fig:Spec-[G:H]=l,k>2}} \\ 
 & & $(\ell-3)k$ & $3\ell-3$ & \quad~~ $\mathcal{P}_{i,j}$, and $\mathcal{S}$ inside $\mathcal{B}^{(1)}$ & \\ 
 & & $(\ell-1)k$ & $3n-3\ell$ & \quad~~ $\mathcal{Q}_{i,j}$, and $\mathcal{U}$ inside $\mathcal{B}^{(2)}$ & \\ 
 & & $-3k$ & $\ell^2-3\ell+2$ & \quad \quad~ $\mathcal{S}$ and $\mathcal{T}$ inside $\mathcal{B}^{(1)}$ & \\ 
 & & $-k$ & $3n\ell-3\ell^2-3n+3\ell$ & \quad~ $\mathcal{R}_{i,j}$, and $\mathcal{U},\mathcal{V}$ inside $\mathcal{B}^{(2)}$ & \\ 
 & & $0$ & $n^2-3n\ell+2\ell^2$ & \quad \quad~ $\mathcal{U}$ and $\mathcal{V}$ inside $\mathcal{B}^{(2)}$ & \\ \hline
\end{tabular}
\caption{{\scriptsize Summary table for the distinct eigenvalues and multiplicities of $\Gamma_{H}(G)$. Here, $n = |G|$, $k = |H|$, and $\ell = [G:H]$.}}
\label{tab:summary-table-eigenvalues}
\end{table}

The proof of Theorem \ref{qsrg-isospectral} proceeds through a sequence of lemmas and propositions. At first we obtain the distinct possible eigenvalues of $\Gamma_H(G)$, for a group $G$ and its normal subgroup $H$ which is followed by the following lemma and a corollary. 

\begin{lemma}\label{eigenvalue-(g,g)}
Let $G$ be a group and $H\le G$. Let $(\rho,V)$ be a representation of $G$. Then the eigenvalues of $\rho(G\setminus H)\coloneqq\sum_{g\in G\setminus H}\rho(g)$ are
$$\left[\underbrace{|G|-|H|,|G|-|H|,\ldots,|G|-|H|}_{\dim\left(V^G\right)~times},\underbrace{-|H|,-|H|,\ldots,-|H|}_{\dim\left(V^H\right)-\dim\left(V^G\right)~times},\underbrace{0,0,\ldots,0}_{\dim V-\dim\left(V^H\right)~times}\right].$$
\end{lemma}
\pf Let 
$$\operatorname{IRR}(G)=\lbrace{(\rho_i,V_i)\mid 1\le i\le t}\rbrace,$$
and 
$$\operatorname{Irr}(G)=\lbrace{\chi_i\mid 1\le i\le t}\rbrace,$$
where $\rho_1$ is the trivial representation and $\chi_1$ is the trivial character of $G$. Then the representation $(\rho,V)$ can be decomposed as
$$\rho\sim\bigoplus_{i=1}^t\rho_i^{\oplus m_i},$$
where $m_i=\langle{\chi_{\rho},\chi_i}\rangle$ is the multiplicity of $\rho_i$ in $\rho$.
Hence,
$$\rho(G\setminus H)\coloneqq\sum_{g\in G\setminus H}\rho(g)\sim\bigoplus_{i=1}^t\left(\sum_{g\in G\setminus H}\rho_i(g)\right)^{\oplus m_i}.$$
Now the eigenvalue of $\rho_1(G\setminus H)$ is $|G|-|H|$ whose multiplicity is $m_1=\langle{\chi_{\rho},\chi_1}\rangle=\dim\left(V^G\right)$. For $2\le i\le t$, $\rho_i(G\setminus H)=-\rho_i(H)$, by Lemma \ref{eigenvalue-projection-operator-irr-lemma} (i). By Lemma \ref{eigenvalue-projection-operator-irr-lemma} (ii), for $2\le i\le t$, the eigenvalues of $\rho_i(G\setminus H)$ are
$$\left[\underbrace{-|H|,-|H|,\ldots,-|H|}_{\dim\left(V_i^H\right)~\text{times}},\underbrace{0,0,\ldots,0}_{\dim V_i-\dim\left(V_i^H\right)~\text{times}}\right].$$
Now since,
$$V=\bigoplus_{i=1}^tV_i^{\oplus m_i},$$
we have,
$$V^H=\bigoplus_{i=1}^t\left(V_i^{H}\right)^{\oplus m_i}.$$
Also, 
$$\dim(V_1^{H})=\frac{1}{|H|}\sum_{h\in H}\chi_1(h)=1.$$ 
Hence we have,
$$\sum_{i=2}^tm_i\dim\left(V_i^H\right)=\dim\left(V^H\right)-\dim\left(V^G\right).$$
Thus the multiplicity of the eigenvalue $-|H|$ in $\rho(G\setminus H)$ becomes $\dim\left(V^H\right)-\dim\left(V^G\right)$ and consequently, the multiplicity of the eigenvalue $0$ is $\dim V-\dim\left(V^H\right)$.\qed

\begin{corollary}\label{corollary-eigenvalue-(g,g)}
Let $G$ be a group and $H\le G$. Let $(\rho,V),(\varphi,W)\in\operatorname{IRR}(G)$. Then the eigenvalues of $(\rho\otimes\varphi)(G\setminus H)\coloneqq\sum_{g\in G\setminus H}(\rho\otimes\varphi)(g)$ are
\begin{enumerate}[label=(\roman*)]
\item $$\left[|G|-|H|,\underbrace{-|H|,-|H|,\ldots,-|H|}_{\dim(V\otimes W)^H-1~times},\underbrace{0,0,\ldots,0}_{\dim (V\otimes W)-\dim(V\otimes W)^H~times}\right],~\text{ if $\rho\sim\varphi^*$}.$$
\item $$\left[\underbrace{-|H|,-|H|,\ldots,-|H|}_{\dim(V\otimes W)^H~times},\underbrace{0,0,\ldots,0}_{\dim (V\otimes W)-\dim(V\otimes W)^H~times}\right],~\text{ if $\rho\nsim\varphi^*$}.$$
\end{enumerate}
\end{corollary}

\begin{proposition}\label{qsrg-dist-eval}
For $H\lhd G$, $\kappa(\Gamma_H(G))\in\lbrace{4,5,6}\rbrace$.
\end{proposition}
\pf Let 
$$\operatorname{IRR}(G)=\lbrace{(\rho_i,V_i)\mid 1\le i\le t}\rbrace,$$
and 
$$\operatorname{Irr}(G)=\lbrace{\chi_i\mid 1\le i\le t}\rbrace,$$
where $\rho_1$ is the trivial representation and $\chi_1$ is the trivial character of $G$. Let $d_i=\deg(\chi_i)$. Since $H$ is a non-trivial, proper subgroup of $G$, $\Gamma_H(G)$ is not strongly regular. By Proposition \ref{simple-eigenval} and Theorem \ref{other-two-eigenval}, $\kappa(\Gamma_H(G))\ge 4$. From \eqref{Omega-S_H}, we have for $1\le i,j\le t$,
$$\Omega_{ij}(S_H)=\frac{\chi_i(G\setminus H)}{d_i}I_{d_id_j}+\frac{\chi_j(G\setminus H)}{d_j}I_{d_id_j}+\Omega_{ij}\left(S_{H}^{(3)}\right).$$
Hence to compute the eigenvalues of $\Omega_{ij}(S_H)$, we need to find the eigenvalues of $$\Omega_{ij}\left(S_{H}^{(3)}\right)=\sum_{g\in G\setminus H}(\rho_i\otimes\rho_j)(g),$$ which can be obtained by Corollary \ref{corollary-eigenvalue-(g,g)}. Let the (multi)-set of the eigenvalues of $\Omega_{ij}(S_H)$ and $\Omega_{ij}\left(S_{H}^{(3)}\right)$ be denoted by $\mathcal{I}_{i,j}$ and $\mathcal{J}_{i,j}$ respectively. Then
\begin{equation}\label{eigenvalue-multiset-eqn}
\mathcal{I}_{i,j}=\left[\underbrace{\frac{\chi_i(G\setminus H)}{d_i},\ldots,\frac{\chi_i(G\setminus H)}{d_i}}_{d_id_j~\text{times}}\right]+\left[\underbrace{\frac{\chi_j(G\setminus H)}{d_j},\ldots,\frac{\chi_j(G\setminus H)}{d_j}}_{d_id_j~\text{times}}\right]+\mathcal{J}_{i,j}.
\end{equation}
Clearly, $\mathcal{J}_{1,1}=\lbrace{|G|-|H|}\rbrace$ and hence $\mathcal{I}_{1,1}=\lbrace{3(|G|-|H|)}\rbrace$. Let 
$$\operatorname{IRR}(G/H)=\lbrace{(\rho_i,V_i)\mid 1\le i\le m}\rbrace,$$
i.e., 
$$\operatorname{IRR}(G)=\left\lbrace{\underbrace{(\rho_1,V_1),(\rho_2,V_2),\ldots,(\rho_m,V_m)}_{\text{trivial~on}~H},\underbrace{(\rho_{m+1},V_{m+1}),\ldots,(\rho_t,V_t)}_{\text{non-trivial~on}~H}}\right\rbrace.$$
In terms of characters, this is equivalent to
$$\operatorname{Irr}(G)=\left\lbrace{\underbrace{\chi_1,\chi_2,\ldots,\chi_m}_{\text{trivial~on}~H},\underbrace{\chi_{m+1},\ldots,\chi_t}_{\text{non-trivial~on}~H}}\right\rbrace.$$
Let $2\le j\le t$. Then $\Omega_{1j}\left(S_{H}^{(3)}\right)=\sum_{g\in G\setminus H}\rho_j(g)$ and hence by Lemma \ref{G-H-char-lemma},
$$\mathcal{J}_{1,j}=\begin{cases}
    \left[\underbrace{-|H|,-|H|,\ldots,-|H|}_{d_j~\text{times}}\right],& \text{if $2\le j\le m$.}\\
    \\
    \left[\underbrace{0,0,\ldots,0}_{d_j~\text{times}}\right],& \text{if $m+1\le j\le t$.}
\end{cases}$$
Thus for $2\le j\le t$, we have from \eqref{eigenvalue-multiset-eqn},
$$\mathcal{I}_{1,j}=\begin{cases}
\left[\underbrace{|G|-3|H|,|G|-3|H|,\ldots,|G|-3|H|}_{d_j~\text{times}}\right],& \text{if $2\le j\le m$.}\\
\\
\left[\underbrace{|G|-|H|,|G|-|H|,\ldots,|G|-|H|}_{d_j~\text{times}}\right],& \text{if $m+1\le j\le t$.}
\end{cases}$$
Similarly from \eqref{eigenvalue-multiset-eqn} we have,
$$\mathcal{I}_{i,1}=\begin{cases}
\left[\underbrace{|G|-3|H|,|G|-3|H|,\ldots,|G|-3|H|}_{d_i~\text{times}}\right],& \text{if $2\le i\le m$}.\\
\\
\left[\underbrace{|G|-|H|,|G|-|H|,\ldots,|G|-|H|}_{d_i~\text{times}}\right],& \text{if $m+1\le i\le t$.}
\end{cases}$$
Now let $2\le i\le m$ and $m+1\le j\le t$. Then by Lemma \ref{closed-under-dual-lemma}, $\rho_i^*\nsim\rho_j$ and thus $\dim(V_i\otimes V_j)^G=0$. Now by Lemma \ref{G-H-char-lemma} (ii),
\begin{align*}
\dim(V_i\otimes V_j)^H&=\frac{1}{|H|}\sum_{h\in H}\chi_i(h)\chi_j
(h)\\
&=\chi_i(1)\left(\frac{1}{|H|}\sum_{h\in H}\chi_j(h)\right)\\
&=0.
\end{align*}
Hence by Corollary \ref{corollary-eigenvalue-(g,g)}, the only eigenvalue of $\Omega_{ij}\left(S_{H}^{(3)}\right)$ is $0$ with multiplicity $d_id_j$. Thus for $2\le i\le m$ and $m+1\le j\le t$, from \eqref{eigenvalue-multiset-eqn} we have,
\begin{align*}
\mathcal{I}_{i,j}
&=\left[\underbrace{-|H|,-|H|,\ldots,-|H|}_{d_id_j~\text{times}}\right].
\end{align*}
Similarly for $2\le j\le m$ and $m+1\le i\le t$, from \eqref{eigenvalue-multiset-eqn} we have,
$$\mathcal{I}_{i,j}=\left[\underbrace{-|H|,-|H|,\ldots,-|H|}_{d_id_j~\text{times}}\right].$$
Now let $2\le i,j\le m$. Then
\begin{align*}
\dim(V_i\otimes V_j)^H&=\frac{1}{|H|}\sum_{h\in H}\chi_i(h)\chi_j
(h)\\
&=\chi_i(1)\chi_j(1)\left(\frac{1}{|H|}\sum_{h\in H}1\right)\\
&=d_id_j.
\end{align*}
Hence by Corollary \ref{corollary-eigenvalue-(g,g)},
$$\mathcal{J}_{i,j}=\begin{cases}
    \left[|G|-|H|,\underbrace{-|H|,-|H|,\ldots,-|H|}_{d_id_j-1~\text{times}}\right],& \text{if $\rho_i^*\sim\rho_j$.}\\
    \\
    \left[\underbrace{-|H|,-|H|,\ldots,-|H|}_{d_id_j~\text{times}}\right],& \text{if $\rho_i^*\nsim\rho_j$.}
\end{cases}$$
which by Lemma \ref{G-H-char-lemma} yields the following from \eqref{eigenvalue-multiset-eqn},
$$\mathcal{I}_{i,j}=\begin{cases}
\left[|G|-3|H|,\underbrace{-3|H|,-3|H|,\ldots,-3|H|}_{d_id_j-1~\text{times}}\right],& \text{if $\rho_i^*\sim\rho_j$.}\\
\\
\left[\underbrace{-3|H|,-3|H|,\ldots,-3|H|}_{d_id_j~\text{times}}\right],& \text{if $\rho_i^*\nsim\rho_j$.}
\end{cases}$$
Finally let $m+1\le i,j\le t$. Then by Corollary \ref{corollary-eigenvalue-(g,g)},
$$\mathcal{J}_{i,j}=\begin{cases}
    \left[|G|-|H|,\underbrace{-|H|,-|H|,\ldots,-|H|}_{\dim(V_i\otimes V_j)^H-1~\text{times}},\underbrace{0,0,\ldots,0}_{d_id_j-\dim(V_i\otimes V_j)^H~\text{times}}\right],& \text{if $\rho_i^*\sim\rho_j$.}\\
    \\
    \left[\underbrace{-|H|,-|H|,\ldots,-|H|}_{\dim(V_i\otimes V_j)^H~\text{times}},\underbrace{0,0,\ldots,0}_{d_id_j-\dim(V_i\otimes V_j)^H~\text{times}}\right],& \text{if $\rho_i^*\nsim\rho_j$.}
\end{cases}$$
which by Lemma \ref{G-H-char-lemma} yields the following from \eqref{eigenvalue-multiset-eqn},
$$\mathcal{I}_{i,j}=\begin{cases}
    \left[|G|-|H|,\underbrace{-|H|,-|H|,\ldots,-|H|}_{\dim(V_i\otimes V_j)^H-1~\text{times}},\underbrace{0,0,\ldots,0}_{d_id_j-\dim(V_i\otimes V_j)^H~\text{times}}\right],& \text{if $\rho_i^*\sim\rho_j$.}\\
    \\
    \left[\underbrace{-|H|,-|H|,\ldots,-|H|}_{\dim(V_i\otimes V_j)^H~\text{times}},\underbrace{0,0,\ldots,0}_{d_id_j-\dim(V_i\otimes V_j)^H~\text{times}}\right],& \text{if $\rho_i^*\nsim\rho_j$.}
\end{cases}$$
Thus the possible eigenvalues of $\Omega_{ij}(S_H)$ are $3(|G|-|H|),|G|-|H|,|G|-3|H|,-|H|,-3|H|,0$. This completes the proof.\qed

\begin{remark}\label{dist-eigenvalue-table-multiset}
The following table describes the eigenvalues of $\Gamma_H(G)$, when $H\lhd G$ in a more comprehensible way. 

$$\begin{array}{|c|cccc|cccc|}
\hline
    & (\rho_1,V_1) & (\rho_2,V_2) & \cdots & (\rho_m,V_m) & (\rho_{m+1},V_{m+1}) & \cdots & (\rho_t,V_t) &\\
\hline
(\rho_1,V_1) & \lbrace{3(|G|-|H|)}\rbrace & \mathcal{P}_{1,2} & \cdots & \mathcal{P}_{1,m} & \mathcal{Q}_{1,m+1} & \cdots & \mathcal{Q}_{1,t} &\\
(\rho_2,V_2) & \mathcal{P}_{2,1} & & & & \mathcal{R}_{2,m+1} & \cdots & \mathcal{R}_{2,t} &\\
\vdots & \vdots & & \mathcal{B}^{(1)} & & \vdots & \ddots & \vdots &\\
(\rho_m,V_m) & \mathcal{P}_{m,1} & & & & \mathcal{R}_{m,m+1} & \cdots & \mathcal{R}_{m,t} &\\
\hline
(\rho_{m+1},V_{m+1}) & \mathcal{Q}_{m+1,1} & \mathcal{R}_{m+1,2} & \cdots & \mathcal{R}_{m+1,m} & & & &\\
\vdots & \vdots & \vdots & \ddots & \vdots & & \mathcal{B}^{(2)} & &\\
(\rho_t,V_t) & \mathcal{Q}_{t,1} & \mathcal{R}_{t,2} & \cdots & \mathcal{R}_{t,m} & & & &\\
\hline
\end{array}$$
where
$$\mathcal{P}_{i,j}=\left[\underbrace{|G|-3|H|,|G|-3|H|,\ldots,|G|-3|H|}_{d_id_j~\text{times}}\right],$$
$$\mathcal{Q}_{i,j}=\left[\underbrace{|G|-|H|,|G|-|H|,\ldots,|G|-|H|}_{d_id_j~\text{times}}\right],$$
$$\mathcal{R}_{i,j}=\left[\underbrace{-|H|,-|H|,\ldots,-|H|}_{d_id_j~\text{times}}\right].$$
And (after suitable rearrangement of the entries),
$$\mathcal{B}^{(1)}=\begin{pmatrix}
    \mathcal{S} & \mathcal{T} & \cdots & \mathcal{T}\\
    \vdots & \vdots & \ddots & \vdots\\
    \mathcal{S} & \mathcal{T} & \cdots & \mathcal{T}
\end{pmatrix},~\mathcal{B}^{(2)}=\begin{pmatrix}
    \mathcal{U} & \mathcal{V} & \cdots & \mathcal{V}\\
    \vdots & \vdots & \ddots & \vdots\\
    \mathcal{U} & \mathcal{V} & \cdots & \mathcal{V}\\
\end{pmatrix}$$
where for $2\le i,j\le m$ with $\rho_i^*\sim\rho_j$,
$$\mathcal{S}=\left[|G|-3|H|,\underbrace{-3|H|,-3|H|,\ldots,-3|H|}_{d_id_j-1~\text{times}}\right],$$
for $2\le i,j\le m$ with $\rho_i^*\nsim\rho_j$,
$$\mathcal{T}=\left[\underbrace{-3|H|,-3|H|,\ldots,-3|H|}_{d_id_j~\text{times}}\right].$$
And for $m+1\le i,j\le t$ with $\rho_i^*\sim\rho_j$,
$$\mathcal{U}=\left [|G|-|H|,\underbrace{-|H|,-|H|,\ldots,-|H|}_{\dim(V_i\otimes V_j)^H-1~\text{times}},\underbrace{0,0,\ldots,0}_{d_id_j-\dim(V_i\otimes V_j)^H~\text{times}}\right],$$
for $m+1\le i,j\le t$ with $\rho_i^*\nsim\rho_j$,
$$\mathcal{V}=\left[\underbrace{-|H|,-|H|,\ldots,-|H|}_{\dim(V_i\otimes V_j)^H~\text{times}},\underbrace{0,0,\ldots,0}_{d_id_j-\dim(V_i\otimes V_j)^H~\text{times}}\right].$$
\end{remark}

\begin{lemma}\label{H-fixed-subsp-dim-sum}
Let $G$ be a group of order $n$ and $H\lhd G$ be such that $[G:H]=\ell$. Let $\operatorname{IRR}(G)=\lbrace{(\rho_i,V_i)\mid 1\le i\le t}\rbrace$ and $\operatorname{IRR}(G/H)=\lbrace{(\rho_i,V_i)\mid 1\le i\le m}\rbrace$. Then
$$\sum_{m+1\le i,j\le t}d_id_j\left(\dim(V_i\otimes V_j)^H\right)=\ell(n-\ell).$$
\end{lemma}
\pf Let 
$$S\coloneqq\sum_{m+1\le i,j\le t}d_id_j\left(\dim(V_i\otimes V_j)^H\right).$$ 
Let $\chi_i$ be the irreducible character of $G$ corresponding to the representation $\rho_i$. Then
\begin{align*}
S&=\sum_{m+1\le i,j\le t}d_id_j\left(\frac{1}{|H|}\sum_{h\in H}\chi_i(h)\chi_j(h)\right)\\
&=\frac{1}{|H|}\sum_{h\in H}\left(\sum_{i=m+1}^td_i\chi_i(h)\right)\left(\sum_{j=m+1}^td_j\chi_j(h)\right).
\end{align*}
Since $H\subseteq\ker(\chi_i)$ for $1\le i\le m$, we have
\begin{align*}
\sum_{i=m+1}^td_i\chi_i(h)&=\sum_{i=1}^td_i\chi_i(h)-\sum_{i=1}^md_i\chi_i(h)\\
&=\chi_L(h)-\sum_{i=1}^md_i^2\\
&=\chi_L(h)-\ell
\end{align*}
where $\chi_L$ is the character of the regular representation. Hence,
\begin{align*}
S&=\frac{1}{|H|}\sum_{h\in H}(\chi_L(h)-\ell)^2\\
&=\frac{1}{|H|}\left[(\chi_L(1)-\ell)^2+\sum_{h\in G, h\ne 1}(\chi_L(h)-\ell)^2\right]\\
&=\frac{1}{|H|}[(n-\ell)^2+(|H|-1)\ell^2].
\end{align*}
The result follows by putting $|H|=n/\ell$.\qed

\medskip

With the help of Proposition \ref{qsrg-dist-eval} and Lemma \ref{H-fixed-subsp-dim-sum}, we compute the complete spectrum of $\Gamma_H(G)$ for a group $G$ and its normal subgroup $H$. For convenience, we shall consistently use $n$ and $k$ to denote $|G|$ and $|H|$, respectively.

\subsection{Spectrum of \texorpdfstring{$\Gamma_H(G)$}{Gamma	extsubscript{H}(G)}, for \texorpdfstring{$[G:H]=2$}{[G:H]=2}}

\medskip

\begin{proposition}\label{index-2-spectrum}
Let $[G:H]=2$. Then $\kappa(\Gamma_H(G))=4$ and
\begin{equation}\label{index-2-spectrum-eqn}
\operatorname{Spec}(\Gamma_H(G))=\begin{pmatrix}
3k & k & -k & 0\\
1 & 3n-6 & 3n-3 & n^2-6n+8
\end{pmatrix}.
\end{equation}
\end{proposition}
\pf Since $[G:H]=2$, the irreducible representations of $G$ can be labelled as,
$$\operatorname{IRR}(G)=\left\lbrace{\underbrace{(\rho_1,V_1),(\rho_2,V_2)}_{\text{trivial~on}~H},\underbrace{(\rho_{3},V_{3}),\ldots,(\rho_t,V_t)}_{\text{non-trivial~on}~H}}\right\rbrace.$$
From Proposition \ref{qsrg-dist-eval} and Remark \ref{dist-eigenvalue-table-multiset}, it follows that $\kappa(\Gamma_H(G))=4$ and $\mathcal{K}(\Gamma_H(G))=\lbrace{3k,k,-k,0}\rbrace$. The eigenvalues of the graph with their corresponding representation-theoretic sources are described in Figure \ref{fig:Spec-[G:H]=2}.

\begin{figure}[htbp]
\centering
$$\begin{array}{|c|cc|ccc|}
\hline
    & (\rho_1,V_1) & (\rho_2,V_2) & (\rho_3,V_3) & \cdots & (\rho_t,V_t)\\
\hline
(\rho_1,V_1) & \lbrace{3k}\rbrace & \lbrace{-k}\rbrace & \mathcal{Q}_{1,3} & \cdots & \mathcal{Q}_{1,t}\\
(\rho_2,V_2) & \lbrace{-k}\rbrace & \lbrace{-k}\rbrace & \mathcal{R}_{2,3} & \cdots & \mathcal{R}_{2,t}\\
\hline
(\rho_3,V_3) & \mathcal{Q}_{3,1} & \mathcal{R}_{3,2} & & &\\
\vdots & \vdots & \vdots & & \mathcal{B}^{(2)} &\\
(\rho_t,V_t) & \mathcal{Q}_{t,1} & \mathcal{R}_{t,2} & & &\\
\hline
\end{array}$$
where
$$\mathcal{Q}_{i,j}=\left[\underbrace{k,k,\ldots,k}_{d_id_j~\text{times}}\right],~\mathcal{R}_{i,j}=\left[\underbrace{-k,-k,\ldots,-k}_{d_id_j~\text{times}}\right].$$
And (after suitable rearrangement of the entries),
$$\mathcal{B}^{(2)}=\begin{pmatrix}
    \mathcal{U} & \mathcal{V} & \cdots & \mathcal{V}\\
    \vdots & \vdots & \ddots & \vdots\\
     \mathcal{U} & \mathcal{V} & \cdots & \mathcal{V}\\
\end{pmatrix}$$
where for $3\le i,j\le t$ with $\rho_i^*\sim\rho_j$,
$$\mathcal{U}=\left[k,\underbrace{-k,-k,\ldots,-k}_{\dim(V_i\otimes V_j)^H-1~\text{times}},\underbrace{0,0,\ldots,0}_{d_id_j-\dim(V_i\otimes V_j)^H~\text{times}}\right],$$
and for $3\le i,j\le t$ with $\rho_i^*\nsim\rho_j$,
$$\mathcal{V}=\left[\underbrace{-k,-k,\ldots,-k}_{\dim(V_i\otimes V_j)^H~\text{times}},\underbrace{0,0,\ldots,0}_{d_id_j-\dim(V_i\otimes V_j)^H~\text{times}}\right].$$
\caption{The eigenvalues of $\Gamma_H(G)$ with multiplicities, when $[G:H]=2$.}
\label{fig:Spec-[G:H]=2}
\end{figure}

\medskip

Thus from \eqref{eqn-2} we have the multiplicities of the eigenvalues as following.
\begin{enumerate}
\item \begin{align*}
m(3k)&=1.
\end{align*}
\item \begin{align*}
m(k)&=\sum_{3\le j\le t}(d_1d_j)^2+\sum_{3\le i\le t}(d_id_1)^2+\sum_{\tiny\begin{array}{cc}
    3\le i,j\le t\\
    \rho_i^*\sim\rho_j 
\end{array}}(d_id_j)\cdot 1\\
&=3\left(\sum_{3\le i\le t}d_i^2\right)\\
&=3(n-2).
\end{align*}
\item \begin{align*}
m(-k)&=(d_1d_2)^2+(d_2d_1)^2+(d_2d_2)^2+\sum_{3\le j\le t}(d_2d_j)^2+\sum_{3\le i\le t}(d_id_2)^2\\
&+\sum_{\tiny\begin{array}{cc}
    3\le i,j\le t\\
    \rho_i^*\sim\rho_j 
\end{array}}d_id_j\left(\dim(V_i\otimes V_j)^H-1\right)+\sum_{\tiny\begin{array}{cc}
    3\le i,j\le t\\
    \rho_i^*\nsim\rho_j 
\end{array}}d_id_j\left(\dim(V_i\otimes V_j)^H\right)\\
&=3+\sum_{3\le j\le t}d_j^2+\sum_{3\le i\le t}d_i^2+\sum_{3\le i,j\le t}d_id_j\left(\dim(V_i\otimes V_j)^H\right)-\sum_{\tiny\begin{array}{cc}
    3\le i,j\le t\\
    \rho_i^*\sim\rho_j 
\end{array}}d_id_j
\end{align*}
\begin{align*}
&=3+(n-2)+2(n-2)~\text{(by Lemma
\ref{H-fixed-subsp-dim-sum})}~~~~~~~~~~~~~~~~~~~~~~~~~~\\
&=3(n-1).
\end{align*}
\item \begin{align*}
m(0)&=n^2-(1+3(n-2)+3(n-1))\\
&=n^2-6n+8.
\end{align*}
\end{enumerate}
This completes the proof.\qed

\begin{example}
Let $G=S_3=\langle{a,b\mid a^3=b^2=1, ab=ba^{-1}}\rangle$ and $H=A_3=\lbrace{1,a,a^2}\rbrace$. Then one can computationally check using SageMath \cite{stein2007sage} that $$\operatorname{Spec}(\Gamma_H(G))=\left[9,[3]^{12},[-3]^{15},[0]^{8}\right],$$
which matches the expression given by \eqref{index-2-spectrum-eqn} in Proposition \ref{index-2-spectrum}.  
\end{example}

\subsection{Spectrum of \texorpdfstring{$\Gamma_H(G)$}{Gamma extsubscript{H}(G)}, for \texorpdfstring{$[G:H]=3$}{[G:H]=3}}

\medskip

\begin{proposition}\label{index-3-spectrum}
Let $H\lhd G$ be such that $[G:H]=3$. Then $\kappa(\Gamma_H(G))=5$ and
\begin{equation}\label{index-3-spectrum-eqn}
\operatorname{Spec}(\Gamma_H(G))=\begin{pmatrix}
6k & 2k & -k & -3k & 0\\
1 & 3n-9 & 6n-18 & 2 & n^2-9n+24
\end{pmatrix}.
\end{equation}
\end{proposition}
\pf Since $[G:H]=3$, the irreducible representations of $G$ can be labelled as,
$$\operatorname{IRR}(G)=\left\lbrace{\underbrace{(\rho_1,V_1),(\rho_2,V_2),(\rho_3,V_3)}_{\text{trivial~on}~H},\underbrace{(\rho_{4},V_{4}),\ldots,(\rho_t,V_t)}_{\text{non-trivial~on}~H}}\right\rbrace.$$
From Proposition \ref{qsrg-dist-eval} and Remark \ref{dist-eigenvalue-table-multiset}, it follows that $\kappa(\Gamma_H(G))=5$ and $\mathcal{K}(\Gamma_H(G))=\lbrace{6k,2k,-k,-3k,0}\rbrace$. The eigenvalues of the graph with their corresponding representation-theoretic sources are described in Figure \ref{fig:Spec-[G:H]=3}.

\begin{figure}[htbp]
\centering
\centering
$$\begin{array}{|c|ccc|ccc|}
\hline
    & (\rho_1,V_1) & (\rho_2,V_2) & (\rho_3,V_3) & (\rho_4,V_4) & \cdots & (\rho_t,V_t)\\
\hline
(\rho_1,V_1) & \lbrace{6k}\rbrace & \lbrace{0}\rbrace & \lbrace{0}\rbrace & \mathcal{Q}_{1,4} & \cdots & \mathcal{Q}_{1,t}\\
(\rho_2,V_2) & \lbrace{0}\rbrace & \lbrace{0}\rbrace & \lbrace{-3k}\rbrace & \mathcal{R}_{2,4} & \cdots & \mathcal{R}_{2,t}\\
(\rho_3,V_3) & \lbrace{0}\rbrace & \lbrace{-3k}\rbrace & \lbrace{0}\rbrace & \mathcal{R}_{3,4} & \cdots & \mathcal{R}_{3,t}\\
\hline
(\rho_4,V_4) & \mathcal{Q}_{4,1} & \mathcal{R}_{4,2} & \mathcal{R}_{4,3} & & &\\
\vdots & \vdots & \vdots & \vdots & & \mathcal{B}^{(2)} &\\
(\rho_t,V_t) & \mathcal{Q}_{t,1} & \mathcal{R}_{t,2} & \mathcal{R}_{t,3} & & &\\
\hline
\end{array}$$
where
$$\mathcal{Q}_{i,j}=\left[\underbrace{2k,2k,\ldots,2k}_{d_id_j~\text{times}}\right],~\mathcal{R}_{i,j}=\left[\underbrace{-k,-k,\ldots,-k}_{d_id_j~\text{times}}\right].$$
And (after suitable rearrangement of the entries),
$$\mathcal{B}^{(2)}=\begin{pmatrix}
    \mathcal{U} & \mathcal{V} & \cdots & \mathcal{V}\\
    \vdots & \vdots & \ddots & \vdots\\
    \mathcal{U} & \mathcal{V} & \cdots & \mathcal{V}\\
\end{pmatrix}$$
where for $4\le i,j\le t$ with $\rho_i^*\sim\rho_j$,
$$\mathcal{U}=\left[2k,\underbrace{-k,-k,\ldots,-k}_{\dim(V_i\otimes V_j)^H-1~\text{times}},\underbrace{0,0,\ldots,0}_{d_id_j-\dim(V_i\otimes V_j)^H~\text{times}}\right],$$
and for $4\le i,j\le t$ with $\rho_i^*\nsim\rho_j$,
$$\mathcal{V}=\left[\underbrace{-k,-k,\ldots,-k}_{\dim(V_i\otimes V_j)^H~\text{times}},\underbrace{0,0,\ldots,0}_{d_id_j-\dim(V_i\otimes V_j)^H~\text{times}}\right].$$
\caption{The eigenvalues of $\Gamma_H(G)$ with multiplicities, when $H\lhd G$ and $[G:H]=3$.}
\label{fig:Spec-[G:H]=3}
\end{figure}

\medskip

Hence from \eqref{eqn-2} the multiplicities of the eigenvalues are as follows.
\begin{enumerate}
\item \begin{align*}
m(6k)&=1.
\end{align*}
\item \begin{align*}
m(2k)&=\sum_{4\le j\le t}(d_1d_j)^2+\sum_{4\le i\le t}(d_id_1)^2+\sum_{\tiny\begin{array}{cc}
    4\le i,j\le t\\
    \rho_i^*\sim\rho_j 
\end{array}}(d_id_j)\cdot 1\\
&=3\left(\sum_{4\le i\le t}d_i^2\right)\\
&=3(n-3).
\end{align*}
\item \begin{align*}
m(-k)&=\sum_{4\le j\le t}(d_2d_j)^2+\sum_{4\le j\le t}(d_3d_j)^2+\sum_{4\le i\le t}(d_id_2)^2+\sum_{4\le i\le t}(d_id_3)^2\\
&+\sum_{\tiny\begin{array}{cc}
    4\le i,j\le t\\
    \rho_i^*\sim\rho_j 
\end{array}}d_id_j\left(\dim(V_i\otimes V_j)^H-1\right)+\sum_{\tiny\begin{array}{cc}
    4\le i,j\le t\\
    \rho_i^*\nsim\rho_j 
\end{array}}d_id_j\left(\dim(V_i\otimes V_j)^H\right)\\
&=4\left(\sum_{4\le i\le t}d_i^2\right)+\sum_{4\le i,j\le t}d_id_j\left(\dim(V_i\otimes V_j)^H\right)-\sum_{\tiny\begin{array}{cc}
    4\le i,j\le t\\
    \rho_i^*\sim\rho_j 
\end{array}}d_id_j\\
&=3\left(\sum_{4\le i\le t}d_i^2\right)+3(n-3)~\text{(by Lemma
\ref{H-fixed-subsp-dim-sum})}\\
&=6(n-3).
\end{align*}
\item \begin{align*}
m(-3k)&=(d_2d_3)^2+(d_3d_2)^2\\
&=2.
\end{align*}
\item \begin{align*}
m(0)&=n^2-(1+3(n-3)+6(n-3)+2)\\
&=n^2-9n+24.
\end{align*}
\end{enumerate}
This completes the proof.\qed

\begin{example}
Let $G=\mathbb{Z}_9$ and $H=\mathbb{Z}_3$. Then one can computationally check using SageMath \cite{stein2007sage} that $$\operatorname{Spec}(\Gamma_H(G))=\left[18,[6]^{18},[-3]^{36},[-9]^{2},[0]^{24}\right],$$ which matches the expression given by \eqref{index-3-spectrum-eqn} in Proposition \ref{index-3-spectrum}.  
\end{example}

\subsection{Spectrum of \texorpdfstring{$\Gamma_H(G)$}{Gamma	extsubscript{H}(G)}, for \texorpdfstring{$[G:H]\ge 4$}{[G:H] >= 4}}

\medskip

Before computing the spectrum of $\Gamma_H(G)$ for index $\ge 4$, we prove a lemma.

\begin{lemma}\label{H-order-2-fixed-subsp-dim}
Let $H\lhd G$ be of order $2$. Let $\operatorname{IRR}(G)=\lbrace{(\rho_i,V_i)\mid 1\le i\le t}\rbrace$ and $\operatorname{IRR}(G/H)=\lbrace{(\rho_i,V_i)\mid 1\le i\le m}\rbrace$. Then for all $m+1\le i,j\le t$,
$$\dim(V_i\otimes V_j)^H=d_id_j.$$
\end{lemma}
\pf Let $H=\lbrace{1,a}\rbrace$. Since $H\lhd G$, by Lemma \ref{G-H-char-lemma} (ii), $\sum_{h\in H}\chi_i(h)=0$ for each $m+1\le i\le t$. Thus we have, $\chi_i(1)+\chi_i(a)=0$, which implies $\chi_i(a)=-\chi_i(1)=-d_i$ for $m+1\le i\le t$. Hence for all $m+1\le i,j\le t$,
\begin{align*}
\dim(V_i\otimes V_j)^H&=\frac{1}{|H|}\sum_{h\in H}\chi_i(h)\chi_j(h)\\
&=\frac{1}{2}(\chi_i(1)\chi_j(1)+\chi_i(a)\chi_j(a))\\
&=\frac{1}{2}(d_id_j+(-d_i)(-d_j))\\
&=d_id_j.
\end{align*}
\qed

\begin{proposition}\label{index-l-spectrum}
Let $H\lhd G$ be such that $[G:H]=\ell\ge 4$.
\begin{enumerate}[label=(\roman*)]
\item If $k=2$, then $\kappa(\Gamma_H(G))=5$ and
\begin{equation}\label{index-l-k=2-spectrum-eqn}
\operatorname{Spec}(\Gamma_H(G))=\begin{pmatrix}
6\ell-6 & 2\ell-6 & 2\ell-2 & -6 & -2\\
1 & 3\ell-3 & 3\ell & \ell^2-3\ell+2 & 3\ell^2-3\ell
\end{pmatrix}.
\end{equation}
\item If $k>2$, then $\kappa(\Gamma_H(G))=6$ and

{\footnotesize
\begin{equation}\label{index-l-k>2-spectrum-eqn}
\operatorname{Spec}(\Gamma_H(G))=\begin{pmatrix}
3(\ell-1)k & (\ell-3)k & (\ell-1)k & -3k & -k & 0\\
1 & 3\ell-3 & 3n-3\ell & \ell^2-3\ell+2 & 3n\ell-3\ell^2-3n+3\ell & n^2-3n\ell+2\ell^2
\end{pmatrix}.
\end{equation}}
\end{enumerate}
\end{proposition}
\pf\\
\textit{Proof of (i):} Since $|H|=k=2$, $|G|=n=2\ell$. Let 
$$\operatorname{IRR}(G)=\lbrace{(\rho_i,V_i)\mid 1\le i\le t}\rbrace,$$
and 
$$\operatorname{IRR}(G/H)=\lbrace{(\rho_i,V_i)\mid 1\le i\le m}\rbrace,$$
i.e.,
$$\operatorname{IRR}(G)=\left\lbrace{\underbrace{(\rho_1,V_1),(\rho_2,V_2),\ldots,(\rho_m,V_m)}_{\text{trivial~on}~H},\underbrace{(\rho_{m+1},V_{m+1}),\ldots,(\rho_t,V_t)}_{\text{non-trivial~on}~H}}\right\rbrace.$$
From Lemma \ref{eigenvalue-(g,g)} and Lemma \ref{H-order-2-fixed-subsp-dim}, it follows that the multiplicity of the eigenvalue $0$ is zero, i.e., $0$ does not appear in the spectrum. From Proposition \ref{qsrg-dist-eval} and Remark \ref{dist-eigenvalue-table-multiset}, it follows that $\kappa(\Gamma_H(G))=5$ and $\mathcal{K}(\Gamma_H(G))=\lbrace{6\ell-6,2\ell-6,2\ell-2,-6,-2}\rbrace$. The eigenvalues of the graph with their corresponding representation-theoretic sources are described in Figure \ref{fig:Spec-[G:H]=l,k=2}.

\begin{figure}[htbp]
\centering
\centering
$$\begin{array}{|c|cccc|ccc|}
\hline
    & (\rho_1,V_1) & (\rho_2,V_2) & \cdots & (\rho_m,V_m) & (\rho_{m+1},V_{m+1}) & \cdots & (\rho_t,V_t)\\
\hline
(\rho_1,V_1) & \lbrace{6(\ell-1)}\rbrace & \mathcal{P}_{1,2} & \cdots & \mathcal{P}_{1,m} & \mathcal{Q}_{1,m+1} & \cdots & \mathcal{Q}_{1,t}\\
(\rho_2,V_2) & \mathcal{P}_{2,1} & & & & \mathcal{R}_{2,m+1} & \cdots & \mathcal{R}_{2,t}\\
\vdots & \vdots & & \mathcal{B}^{(1)} & & \vdots & \ddots & \vdots\\
(\rho_m,V_m) & \mathcal{P}_{m,1} & & & & \mathcal{R}_{m,m+1} & \cdots & \mathcal{R}_{m,t}\\
\hline
(\rho_{m+1},V_{m+1}) & \mathcal{Q}_{m+1,1} & \mathcal{R}_{m+1,2} & \cdots & \mathcal{R}_{m+1,m} & & &\\
\vdots & \vdots & \vdots & \ddots & \vdots & & \mathcal{B}^{(2)} &\\
(\rho_t,V_t) & \mathcal{Q}_{t,1} & \mathcal{R}_{t,2} & \cdots & \mathcal{R}_{t,m} & & &\\
\hline
\end{array}$$
where
$$\mathcal{P}_{i,j}=\left[\underbrace{2\ell-6,2\ell-6,\ldots,2\ell-6}_{d_id_j~\text{times}}\right],~\mathcal{Q}_{i,j}=\left[\underbrace{2\ell-2,2\ell-2,\ldots,2\ell-2}_{d_id_j~\text{times}}\right],$$
$$\mathcal{R}_{i,j}=\left[\underbrace{-2,-2,\ldots,-2}_{d_id_j~\text{times}}\right].$$
And (after suitable rearrangement of the entries),
$$\mathcal{B}^{(1)}=\begin{pmatrix}
    \mathcal{S} & \mathcal{T} & \cdots & \mathcal{T}\\
    \vdots & \vdots & \ddots & \vdots\\
    \mathcal{S} & \mathcal{T} & \cdots & \mathcal{T}
\end{pmatrix},~\mathcal{B}^{(2)}=\begin{pmatrix}
    \mathcal{U} & \mathcal{V} & \cdots & \mathcal{V}\\
    \vdots & \vdots & \ddots & \vdots\\
    \mathcal{U} & \mathcal{V} & \cdots & \mathcal{V}\\
\end{pmatrix}$$
where for $2\le i,j\le m$,
$$\mathcal{S}=\left[2\ell-6,\underbrace{-6,-6,\ldots,-6}_{d_id_j-1~\text{times}}\right],~\mathcal{T}=\left[\underbrace{-6,-6,\ldots,-6}_{d_id_j~\text{times}}\right].$$
And for $m+1\le i,j\le t$,
$$\mathcal{U}=\left[2\ell-2,\underbrace{-2,-2,\ldots,-2}_{d_id_j-1~\text{times}}\right],~\mathcal{V}=\left[\underbrace{-2,-2,\ldots,-2}_{d_id_j~\text{times}}\right].$$
\caption{The eigenvalues of $\Gamma_H(G)$ with multiplicities, when $H\lhd G$, $[G:H]\ge 4$ and $|H|=2$.}
\label{fig:Spec-[G:H]=l,k=2}
\end{figure}

\medskip

Hence from \eqref{eqn-2} the multiplicities of the eigenvalues are as follows.
\begin{enumerate}
\item \begin{align*}
m(6\ell-6)=1.
\end{align*}
\item \begin{align*}
m(2\ell-6)&=\sum_{2\le j\le m}(d_1d_j)^2+\sum_{2\le i\le m}(d_id_1)^2+\sum_{\tiny\begin{array}{cc}
     2\le i,j\le m\\
     \rho_i^*\sim\rho_j
\end{array}}(d_id_j)\cdot 1\\
&=3\left(\sum_{2\le i\le m}d_i^2\right)\\
&=3(\ell-1).
\end{align*}
\item \begin{align*}
m(2\ell-2)&=\sum_{m+1\le j\le t}(d_1d_j)^2+\sum_{m+1\le i\le t}(d_id_1)^2+\sum_{\tiny\begin{array}{cc}
     m+1\le i,j\le t\\
     \rho_i^*\sim\rho_j
\end{array}}(d_id_j)\cdot 1\\
&=3\left(\sum_{m+1\le i\le t}d_i^2\right)\\
&=3\left(\sum_{1\le i\le t}d_i^2-\sum_{1\le i\le m}d_i^2\right)\\
&=3(n-\ell)\\
&=3\ell.
\end{align*}
\item \begin{align*}
m(-6)&=\sum_{\tiny\begin{array}{cc}
    2\le i,j\le m\\
     \rho_i^*\sim\rho_j 
\end{array}}d_id_j(d_id_j-1)+\sum_{\tiny\begin{array}{cc}
    2\le i,j\le m\\
     \rho_i^*\nsim\rho_j 
\end{array}}(d_id_j)(d_id_j)\\
&=\sum_{2\le i,j\le m}d_i^2d_j^2-\sum_{\tiny\begin{array}{cc}
    2\le i,j\le m\\
     \rho_i^*\sim\rho_j 
\end{array}}d_id_j\\
&=\sum_{1\le i,j\le m}d_i^2d_j^2-\sum_{2\le j\le m}d_1^2d_j^2-\sum_{2\le i\le m}d_i^2d_1^2-d_1^2d_1^2-\sum_{2\le i\le m}d_i^2\\
&=\ell^2-3\left(\sum_{2\le i\le m}d_i^2\right)-1\\
&=\ell^2-3\ell+2.  
\end{align*}
\item \begin{align*}
m(-2)&=n^2-(1+3(\ell-1)+3\ell+\ell^2-3\ell+2)\\
&=3\ell^2-3\ell.
\end{align*}
\end{enumerate}
\textit{Proof of (ii):} From Proposition \ref{qsrg-dist-eval} and Remark \ref{dist-eigenvalue-table-multiset}, it follows that $\kappa(\Gamma_H(G))=6$ and $\mathcal{K}(\Gamma_H(G))=\lbrace{3(\ell-1)k,(\ell-3)k,(\ell-1)k,-3k,-k,0}\rbrace$. The eigenvalues of the graph with their corresponding representation-theoretic sources are described in Figure \ref{fig:Spec-[G:H]=l,k>2}.

\begin{figure}[htbp]
\centering
\centering
$$\begin{array}{|c|cccc|cccc|}
\hline
    & (\rho_1,V_1) & (\rho_2,V_2) & \cdots & (\rho_m,V_m) & (\rho_{m+1},V_{m+1}) & \cdots & (\rho_t,V_t) &\\
\hline
(\rho_1,V_1) & \lbrace{3(\ell-1)k}\rbrace & \mathcal{P}_{1,2} & \cdots & \mathcal{P}_{1,m} & \mathcal{Q}_{1,m+1} & \cdots & \mathcal{Q}_{1,t} &\\
(\rho_2,V_2) & \mathcal{P}_{2,1} & & & & \mathcal{R}_{2,m+1} & \cdots & \mathcal{R}_{2,t} &\\
\vdots & \vdots & & \mathcal{B}^{(1)} & & \vdots & \ddots & \vdots &\\
(\rho_m,V_m) & \mathcal{P}_{m,1} & & & & \mathcal{R}_{m,m+1} & \cdots & \mathcal{R}_{m,t} &\\
\hline
(\rho_{m+1},V_{m+1}) & \mathcal{Q}_{m+1,1} & \mathcal{R}_{m+1,2} & \cdots & \mathcal{R}_{m+1,m} & & & &\\
\vdots & \vdots & \vdots & \ddots & \vdots & & \mathcal{B}^{(2)} & &\\
(\rho_t,V_t) & \mathcal{Q}_{t,1} & \mathcal{R}_{t,2} & \cdots & \mathcal{R}_{t,m} & & & &\\
\hline
\end{array}$$
where
$$\mathcal{P}_{i,j}=\left[\underbrace{(\ell-3)k,(\ell-3)k,\ldots,(\ell-3)k}_{d_id_j~\text{times}}\right],~\mathcal{Q}_{i,j}=\left[\underbrace{(\ell-1)k,(\ell-1)k,\ldots,(\ell-1)k}_{d_id_j~\text{times}}\right],$$
$$\mathcal{R}_{i,j}=\left[\underbrace{-k,-k,\ldots,-k}_{d_id_j~\text{times}}\right].$$
And (after suitable rearrangement of the entries),
$$\mathcal{B}^{(1)}=\begin{pmatrix}
    \mathcal{S} & \mathcal{T} & \cdots & \mathcal{T}\\
    \vdots & \vdots & \ddots & \vdots\\
    \mathcal{S} & \mathcal{T} & \cdots & \mathcal{T}
\end{pmatrix},~\mathcal{B}^{(2)}=\begin{pmatrix}
    \mathcal{U} & \mathcal{V} & \cdots & \mathcal{V}\\
    \vdots & \vdots & \ddots & \vdots\\
    \mathcal{U} & \mathcal{V} & \cdots & \mathcal{V}\\
\end{pmatrix}$$
where for $2\le i,j\le m$,
$$\mathcal{S}=\left[(\ell-3)k,\underbrace{-3k,-3k,\ldots,-3k}_{d_id_j-1~\text{times}}\right],~\mathcal{T}=\left[\underbrace{-3k,-3k,\ldots,-3k}_{d_id_j~\text{times}}\right].$$
And for $m+1\le i,j\le t$,
$$\mathcal{U}=\left[(\ell-1)k,\underbrace{-k,-k,\ldots,-k}_{\dim(V_i\otimes V_j)^H-1~\text{times}},\underbrace{0,0,\ldots,0}_{d_id_j-\dim(V_i\otimes V_j)^H~\text{times}}\right],$$
$$\mathcal{V}=\left[\underbrace{-k,-k,\ldots,-k}_{\dim(V_i\otimes V_j)^H~\text{times}},\underbrace{0,0,\ldots,0}_{d_id_j-\dim(V_i\otimes V_j)^H~\text{times}}\right].$$
\caption{The eigenvalues of $\Gamma_H(G)$ with multiplicities, when $H\lhd G$, $[G:H]\ge 4$ and $|H|>2$.}
\label{fig:Spec-[G:H]=l,k>2}
\end{figure}

\medskip

Hence from \eqref{eqn-2} the multiplicities of the eigenvalues are as follows.
\begin{enumerate}
\item \begin{align*}
m(3(\ell-1)k)=1.
\end{align*}
\item \begin{align*}
m((\ell-3)k)&=\sum_{2\le j\le m}(d_1d_j)^2+\sum_{2\le i\le m}(d_id_1)^2+\sum_{\tiny\begin{array}{cc}
     2\le i,j\le m\\
     \rho_i^*\sim\rho_j
\end{array}}(d_id_j)\cdot 1\\
&=3\left(\sum_{2\le i\le m}d_i^2\right)\\
&=3(\ell-1).
\end{align*}
\item \begin{align*}
m((\ell-1)k)&=\sum_{m+1\le j\le t}(d_1d_j)^2+\sum_{m+1\le i\le t}(d_id_1)^2+\sum_{\tiny\begin{array}{cc}
     m+1\le i,j\le t\\
     \rho_i^*\sim\rho_j
\end{array}}(d_id_j)\cdot 1\\
&=3\left(\sum_{m+1\le i\le t}d_i^2\right)\\
&=3\left(\sum_{1\le i\le t}d_i^2-\sum_{1\le i\le m}d_i^2\right)\\
&=3(n-\ell).
\end{align*}
\item \begin{align*}
m(-3k)&=\sum_{\tiny\begin{array}{cc}
    2\le i,j\le m\\
     \rho_i^*\sim\rho_j 
\end{array}}d_id_j(d_id_j-1)+\sum_{\tiny\begin{array}{cc}
    2\le i,j\le m\\
     \rho_i^*\nsim\rho_j 
\end{array}}(d_id_j)(d_id_j)\\
&=\sum_{2\le i,j\le m}d_i^2d_j^2-\sum_{\tiny\begin{array}{cc}
    2\le i,j\le m\\
     \rho_i^*\sim\rho_j 
\end{array}}d_id_j\\
&=\sum_{1\le i,j\le m}d_i^2d_j^2-\sum_{2\le j\le m}d_1^2d_j^2-\sum_{2\le i\le m}d_i^2d_1^2-d_1^2d_1^2-\sum_{2\le i\le m}d_i^2\\
&=\ell^2-3\left(\sum_{2\le i\le m}d_i^2\right)-1\\
&=\ell^2-3(\ell-1)-1\\
&=\ell^2-3\ell+2.
\end{align*}
\item \begin{align*}
m(-k)&=\sum_{\tiny\begin{array}{cc}
    2\le i\le m\\
    m+1\le j\le t 
\end{array}}d_i^2d_j^2+\sum_{\tiny\begin{array}{cc}
    2\le j\le m\\
    m+1\le i\le t 
\end{array}}d_i^2d_j^2+\sum_{\tiny\begin{array}{cc}
     m+1\le i,j\le t\\
     \rho_i^*\sim\rho_j
\end{array}}d_id_j\left(\dim(V_i\otimes V_j)^H-1\right)\\
&+\sum_{\tiny\begin{array}{cc}
     m+1\le i,j\le t\\
     \rho_i^*\nsim\rho_j
\end{array}}d_id_j\left(\dim(V_i\otimes V_j)^H\right)\\
&=2\left(\sum_{\tiny\begin{array}{cc}
    2\le i\le m\\
    m+1\le j\le t 
\end{array}}d_i^2d_j^2\right)+\sum_{m+1\le i,j\le t}d_id_j\left(\dim(V_i\otimes V_j)^H\right)-\sum_{\tiny\begin{array}{cc}
     m+1\le i,j\le t\\
     \rho_i^*\sim\rho_j
\end{array}}d_id_j\\
&=2\left(\sum_{2\le i\le m}d_i^2\right)\left(\sum_{m+1\le j\le t}d_j^2\right)+\sum_{m+1\le i,j\le t}d_id_j\left(\dim(V_i\otimes V_j)^H\right)-\sum_{m+1\le i\le t}d_i^2\\
&=2(\ell-1)(n-\ell)+\ell(n-\ell)-(n-\ell)~\text{(by Lemma \ref{H-fixed-subsp-dim-sum})}\\
&=3n\ell-3\ell^2-3n+3\ell.
\end{align*}
\item \begin{align*}
m(0)&=n^2-(1+3(\ell-1)+3(n-\ell)+\ell^2-3\ell+2+(3n\ell-3\ell^2-3n+3\ell))\\
&=n^2-3n\ell+2\ell^2.
\end{align*}
\end{enumerate}
This completes the proof.\qed

\begin{example}
Let $G=D_4=\langle{a,b\mid a^4=b^2=1, ab=ba^{-1}}\rangle$ and $H=Z(D_4)=\lbrace{1,a^2}\rbrace$. Then one can computationally check using SageMath \cite{stein2007sage} that $$\operatorname{Spec}(\Gamma_H(G))=\left[18,[2]^{9},[6]^{12},[-6]^{6},[-2]^{36}\right],$$ which matches the expression given by \eqref{index-l-k=2-spectrum-eqn} in Proposition \ref{index-l-spectrum}.  
\end{example}

\begin{example}
Let $G=\mathbb{Z}_{18}$ and $H=\mathbb{Z}_3$. Then one can computationally check using SageMath \cite{stein2007sage} that $$\operatorname{Spec}(\Gamma_H(G))=\left[45,[9]^{15},[15]^{36},[-9]^{20},[-3]^{180},[0]^{72}\right],$$ which matches the expression given by \eqref{index-l-k>2-spectrum-eqn} in Proposition \ref{index-l-spectrum}.  
\end{example}

\medskip

Propositions \ref{index-2-spectrum},\ref{index-3-spectrum} and \ref{index-l-spectrum} yield the following theorem.

\begin{theorem}\label{qsrg-isospectral}
Let $G_1$ and $G_2$ be two groups of same order. Let $H_1\lhd G_1$ and $H_2\lhd G_2$ be such that $|H_1|=|H_2|$. Then the graphs $\Gamma_{H_1}(G_1)$ and $\Gamma_{H_2}(G_2)$ are isospectral.\qed
\end{theorem}

\begin{remark}
In Theorem \ref{qsrg-isospectral}, the normality condition of the subgroups is necessary. For example, take $G_1=S_3=\langle{a,b\mid a^3=b^2=1, ab=ba^{-1}}\rangle$, $H_1=\lbrace{1,b}\rbrace$ and $G_2=\mathbb{Z}_6$, $H_2=\mathbb{Z}_2$. Using SageMath \cite{stein2007sage}, one can show that
$$\operatorname{Spec}(\Gamma_{H_1}(G_1))=\begin{pmatrix}
    12 & 0 & 4 & -2 & -1+\sqrt{13} & -1-\sqrt{13}\\
    1 & 4 & 7 & 16 & 4 & 4
\end{pmatrix},$$
and
$$\operatorname{Spec}(\Gamma_{H_2}(G_2))=\begin{pmatrix}
    12 & 0 & 4 & -2 & -6\\
    1 & 6 & 9 & 18 & 2
\end{pmatrix}.$$
\end{remark}

\begin{remark} 
It is not necessarily true that $\Gamma_{H_1}(G_1)\cong\Gamma_{H_2}(G_2)$, even if the conditions in Theorem \ref{qsrg-isospectral} are met. For example, take $G_1=\mathbb{Z}_8$, $H_1=\mathbb{Z}_4$ and $G_2=\mathbb{Z}_4\times\mathbb{Z}_2$, $H_2=\lbrace{(0,0),(0,1),(2,0),(2,1)}\rbrace\cong\mathbb{Z}_2\times\mathbb{Z}_2$. Then using SageMath \cite{stein2007sage}, one can show that $\Gamma_{H_1}(G_1)\ncong\Gamma_{H_2}(G_2)$.
\end{remark}

\section{Conclusion and Open Issues}

In this paper, we identified an integral family of quasi-strongly regular graphs $\Gamma_H(G)$ of grade greater than $2$. We provided a sufficient condition for the integrality of the graph $\Gamma_H(G)$. We also computed the complete spectrum of $\Gamma_H(G)$ when $H\lhd G$. This work opens up several promising avenues for future research. We conclude by outlining some key open problems:
\begin{itemize}
\item Computational evidence using SageMath suggests that the normality of $H$ is not only sufficient but also necessary for the integrality of $\Gamma_H(G)$. However, a formal proof of this statement remains an intriguing challenge.
\item While we have computed the full spectrum for normal subgroups, the spectral properties of $\Gamma_H(G)$, apart from the ones obtained in Theorem \ref{other-two-eigenval}, are largely unexplored. A general characterization of the spectrum in this case would be another interesting area of further investigation.
\end{itemize}

\section*{Acknowledgement}

The first author is supported by the funding of UGC [NTA Ref. No. 211610129182], Govt. of India. The second author is supported by the post-doctoral fellowship of Indian Institute of Technology Bombay, India. The third author acknowledges the funding of DST-FIST Sanction no. $SR/FST/MS-I/2019/41$, Govt. of India.

The authors are grateful to the anonymous referee for his/her insightful comments. The authors would also like to thank Prof. Amritanshu Prasad of Institute of Mathematical Sciences, India, for some fruitful discussion regarding Theorem \ref{H-normal-graph-integral}.

The authors used generative AI tools during the preparation of the manuscript for language editing and for assistance in checking and reformulating parts of the exposition. All mathematical statements and proofs were independently verified by the authors.

\bibliographystyle{abbrv}
\bibliography{ref}

\section*{Statements and Declarations}
The authors have no relevant financial or non-financial interests to disclose. The authors have no competing interests to declare that are relevant to the content of this article. The authors also declare that there is no associated data.

All authors contributed to the study conception and design. Material preparation and analysis were performed by Sauvik Poddar, Sucharita Biswas and Angsuman Das. The first draft of the manuscript was written by Sauvik Poddar, and all the authors commented on previous versions of the manuscript. All authors read and approved the final manuscript.

\end{document}